\documentclass[onefignum,onetabnum]{siamart190516}

\usepackage{lipsum}
\usepackage{graphicx}
\usepackage{epstopdf}
\usepackage[english]{babel}
\usepackage{amsmath,amssymb,latexsym,amsfonts}
\usepackage{bm}
\usepackage{bbm}
\usepackage{soul}
\usepackage{subcaption}
\usepackage{mathtools}
\usepackage{mathrsfs}
\usepackage{xcolor}
\usepackage{algorithm}
\usepackage{algpseudocode}

\ifpdf%
  \DeclareGraphicsExtensions{.eps,.pdf,.png,.jpg}
\else
  \DeclareGraphicsExtensions{.eps}
\fi
\usepackage{amsopn}

\newcommand{\argmin}{\mathop{\mathrm{argmin}}}

\usepackage{booktabs}

\newlength{\tabcolsepold}
\begin{document}
\headers{A Posteriori Error Estimates for Solving Graph Laplacians}{X.~Hu, K.~Wu, and L.~T.~Zikatanov}
\title{A Posteriori Error Estimates for Solving Graph Laplacians \thanks{The work of Hu and Wu was partially supported by the National Science Foundation under grant DMS-1812503 and CCF-1934553. The work of Zikatanov was supported in part by NSF DMS-1720114 and DMS-1819157}}

\author{Xiaozhe Hu \thanks{Department of Mathematics, Tufts University, Medford, MA 02155, USA \, (\email{Xiaozhe.Hu@tufts.edu})}
	\and Kaiyi Wu \thanks{Department of Mathematics, Tufts University, Medford, MA 02155, USA\,(\email{kaiyi.wu@tufts.edu})}
	\and Ludmil T. Zikatanov \thanks{Department of Mathematics, The Pennsylvania State University, University Park, PA 16802, USA\,(\email{ludmil@psu.edu}).}
}

\maketitle
 
\begin{abstract}
 In this paper, we study a posteriori error estimators which aid multilevel iterative solvers for linear systems with graph Laplacians. In earlier works such estimates were computed by solving global optimization problems, which could be computationally expensive. We propose a novel strategy to compute these estimates by constructing a Helmholtz decomposition on the graph based on a spanning tree and the corresponding cycle space. To compute the error estimator, we solve efficiently the linear system on the spanning tree, and then we solve approximately a least-squares problem on the cycle space. As we show, such an estimator has a nearly-linear computational complexity for sparse graphs under certain assumptions.  Numerical experiments are presented to demonstrate the efficacy of the proposed method.
\end{abstract}


\section{Introduction}\label{sec:intro}
Graphs are frequently employed to model networks in social science, energy, and biological applications~\cite{Borgatti09_SocialNetworks, Bullmore09_BioNetworks,Gutman72_EnergyNetworks}. In many cases, these applications require the solution of large-scale linear systems of equations given by the graph Laplacian matrix~\cite{Hu2017_GraphLaplacian,Hu2018_GraphLaplacian,Merris1994143,2012Gerta_ElectricalNetwork,Weinberger2006_GraphLaplacian}. Such linear systems have also been the key to develop link-based ranking algorithms for web searching queries~\cite{liu2011noise, page1999pagerank,  wang2012multimodal} or recommendation systems~\cite{ fouss2007random,ma2016diffusion}. Exploiting the properties of graph Laplacians to achieve dimension reduction in high dimensional data representation, researchers have produced fruitful results in image classification~\cite{gomez2008semisupervised,yu2012adaptive}, representation learning~\cite{hamilton2017representation}, and clustering~\cite{belkin2001laplacian, nie2016constrained}. In the numerical solutions of partial differential equations (PDEs), the stiffness matrices arising from the finite-element or finite-difference method also take the form of graph Laplacian as discussed in \cite{Xu2018_Agg}.  Therefore, it is important to develop efficient and robust methods for solving graph Laplacian systems.

To solve large-scale graph Laplacian linear systems, direct methods suffer from their expensive computational costs~\cite{Zikatanov2003_AMGgraph}. Iterative methods, such as the algebraic multigrid (AMG) methods originated in~\cite{1stAMG}, are often applied to solve the linear systems (see also~\cite{Zikatanov17_AMGSurvey} and the references therein for a recent survey on the AMG methods). In practice, the AMG method achieves optimal computational complexity for many applications, including solving the linear systems with weighted graph Laplacians~\cite{Bolten.M;Friedhoff.S;Frommer.A;Heming.M;Kahl.K.2011a,Zikatanov2013_AMGgraph,DAmbra.P;Vassilevski.P2015a,Zikatanov2003_AMGgraph,koutis2011combinatorial, Livne.O;Brandt.A.2012a, Napov.A;Notay.Y2016a}.

A typical multigrid method traverses a hierarchy of spaces (grids) of different dimensions (grid sizes) and solves the corresponding linear system at different resolutions. As is well known, an efficient AMG method should damp the algebraic high-frequency error using relaxations/smoothers (e.g. common iterative solvers like Jacobi and Gauss Seidel methods) and eliminate the algebraically smooth (low-frequency) error via a coarse space (grid) correction~\cite{Briggs.W;Henson.V;McCormick.S.2000a}. The latter requires the ``smooth'' error obtained on the fine grids to be accurately approximated on the coarse spaces. Many different coarsening strategies have been developed based on good estimations of the error, for example, the classical AMG~\cite{1stAMG,2006Falgout_Adaptive}, the smoothed aggregation AMG~\cite{Brezina.M;Falgout.R;Maclachlan.S;Manteuffel.T;Mccormick.S;Ruge.J.2005a,Brezina.M;Falgout.R;MacLachlan.S;Manteuffel.T;McCormick.S;Ruge.J.2004a,Falgout2010_ASAM}, bootstrap AMG \cite{BS_2015,Brandt2011_BootstrapAMG}, and the unsmoothed aggregation AMG~\cite{Zikatanov2013_AMGgraph2,Zikatanov2013_AMGgraph, Zikatanov2003_AMGgraph,Livne.O;Brandt.A.2012a,Napov.A;Notay.Y.2012a,Notay.Y.2010b,Zikatanov2015_AMGgraph}.  Thus, an efficient, reliable, and computable estimate of the error a posteriori, during AMG iterations, is needed for developing robust AMG methods.

Generally, the a posteriori estimators provide a computable estimation for the true error locally. Our approach borrows several ideas from the finite-element (FE) literature (equilibrated error estimators~\cite{Ainsworth1997_APEEFEM,Rheinboldt1978_APEEFEM, Kelly1983_APEEFEM,Verfurth1994_APEEFEM} and functional a posteriori error estimators~\cite{ APEE2,APEE1, APEE4,APEE3}). In \cite{Xu2018_Agg}, the authors derived a posteriori error estimator for solving graph Laplacian linear systems for the first time based on the functional a posteriori error estimation framework. Such a technique was used to predict the error of approximation from coarse grids for multilevel unsmoothed aggregation AMG and the estimator is computed by solving a perturbed global optimization problem (discussed in more details in~\cref{sec:notations}). Such an approach provides an accurate error estimator. However, it could be computationally expensive, which affects the effectiveness of the resulting adaptive AMG method. In this work, we propose a novel a posteriori error estimator and an efficient algorithm to reduce the computational cost, which could further be used to construct the efficient multilevel hierarchy for adaptive AMG. Roughly speaking, this is achieved by taking advantage of the Helmholtz decomposition~\cite{koenigsberger1902hermann} on the graph computationally, which splits the error into a divergence-free component and a curl-free component. The rationale of the proposed algorithm for computing the approximation to the true error has two main steps:
\begin{enumerate}
	\item
	Solve a linear system on a spanning tree  (defined in~\cref{sec:notations}) of the graph to get the curl-free component, or equivalently, the $\operatorname{grad}$ component of the error. 
	\item
          Approximately solve a minimization problem in the cycle space to obtain the div-free component of the Helmholtz decomposition of the error.
\end{enumerate}
The first step can be done in linear time with Gaussian elimination using lexicographical ordering as shown in~\cite{1976RoseD_TarjanR_LuekerG-aa,2014VassilevskiP_ZikatanovL-aa}.  Solving exactly the minimization in the cycle space of the graph is computationally expensive as it is equivalent to solving a constrained minimization problem, which is as difficult as the original linear system (often even more difficult). Our algorithm solves the minimization approximately by applying several steps of a relaxation scheme, the one level Schwarz method~\cite{2020Hu_SubspaceCorrection,2002Xu_SubspaceCorrection}. This crucial improvement reduces the computational cost and gives an accurate a posteriori estimates of the true error in nearly optimal time for sparse graphs, which is further verified by our numerical experiments.  Clearly, such an error estimator can be incorporated to construct multilevel hierarchies since it provides accurate estimate of the true error, which is important for constructing coarse levels adaptively. The corresponding adaptive coarsening scheme will preserve the smooth error accurately on the coarse levels and ensure the robustness of the resulting adaptive AMG method.   

The rest of this paper is organized as follows. In~\cref{sec:notations} we review backgrounds on graphs and graph Laplacians, along with some previous results in \cite{Xu2018_Agg}. The main algorithm to compute a posterior error estimates is stated in~\cref{sec: main}. We present and analyze some numerical experiments in~\cref{sec:num}. Finally, in~\cref{sec:conc} we summarize the main contribution and list some future work.

\section{Preliminaries}\label{sec:notations}
In this section, we define the necessary notations and recall some fundamental results for the computation of an a posteriori error estimator for solving graph Laplacians as presented in~\cite{Xu2018_Agg}.

\subsection{Graphs and Graph Laplacians}
Consider an undirected weighted graph $\mathcal{G} = (\mathcal{V}, \mathcal{E}, \omega)$, where $\mathcal{V} = \big\{1,2,\cdots, n\big\}$ is the vertex set, $\mathcal{E} = \big\{ \big\{ i,j  \big\},\; i,j \in \mathcal{V}, i > j  \big\}$ is the edge set, and $ \omega = \{\omega_e\}_{e\in \mathcal{E}}$ is the set of edge weights. Here, the weights are assumed to be positive, i.e., $\omega_e > 0$, and  we take all the edge weights to be $1$ for unweighted graphs. We only consider undirected graphs here and that is why we have fixed $i>j$ for every
$e  = \big\{i,j\big\} \in \mathcal{E}$. Thus $\{1,2\}$ cannot be an edge in our graph, while $\{2,1\}$ could be in $\mathcal{E}$. 

Denote $n =| \mathcal{V}| $ and $m = |\mathcal{E}|$.  Let $ \mathscr{V}= \mathbb{R}^n$ and $\mathscr{W}= \mathbb{R}^m$ be the vertex space and edge space, respectively. The inner product on vertex space and edge space are defined as:
\begin{equation*}
\begin{split}
(\bm{u},\bm{v}) = \bm{v}^\intercal \bm{u}, \quad \forall \ \bm{u},\bm{v} \in \mathscr{V}, \\
(\bm{ \tau}, \bm{\phi}) = \bm{\phi}^\intercal \bm{\tau}, \quad \forall \ \bm{ \tau}, \bm{\phi} \in \mathscr{W}.
\end{split}
\end{equation*}
 The weighted graph Laplacian matrix $L \in \mathbb{R}^{n\times n}$ can be defined via the bilinear form:
\begin{equation*} \label{eqn: graph_Laplacian_L}
(L\bm{u},\bm{v})  := \bm{v}^\intercal L\bm{u} = \sum_{ e= \{  i,j \} \in \mathcal{E}} \omega_e(\bm{u}_i-\bm{u}_j)(\bm{v}_i-\bm{v}_j), \quad \forall \; \bm{u},\bm{v} \in \mathscr{V}.
\end{equation*}
Associated with the graph is the discrete gradient operator (or edge-node incidence matrix) $G \in \mathbb{R}^{m \times n}: \mathscr{V} \rightarrow \mathscr{W}$ and the edge weight matrix $D \in \mathbb{R}^{m \times m}: \mathscr{W} \rightarrow \mathscr{W}$. They are defined as the following: for each edge $e  = \big\{i,j\big\} \in \mathcal{E}$,
\begin{equation} \label{eqn: operators}
\begin{split}
(G\bm{v})_e &= \bm{v}_{i}-\bm{v}_{j}, \quad \forall \ \bm{v} \in \mathscr{V}, \\
(D \bm{\tau})_e &= w_e \bm{\tau}_e, \quad \forall \ \bm{\tau} \in  \mathscr{W}.
\end{split}
\end{equation}
The adjoint of $G$, denoted by $G^\intercal: \mathscr{W} \rightarrow \mathscr{V}$, is the discrete divergence operator (or node-edge incidence matrix) on the graph,
\begin{equation} \label{eqn: operator_divergence}
(G\bm{u}, \bm{\tau}) = (\bm{u}, G^\intercal \bm{\tau}), \quad \forall \ \bm{u} \in \mathscr{V}, \; \forall \ \bm{\tau} \in \mathscr{W}.
\end{equation}
By direct computation, the following identity holds true,
\begin{equation*}
(L\bm{u},\bm{v}) = (DG\bm{u}, G\bm{v}).
\end{equation*}
Thus, we can write $L := G^\intercal D G$.  Based on this definition of the graph Laplacian $L$, it is straightforward to verify that, 
\begin{equation*}
\|DG\bm{u}-DG\bm{v}\|_{D^{-1}}^2 = \|\bm{u}-\bm{v}\|^2_L, \quad \forall \ \bm{u},\bm{v} \in \mathscr{V},
\end{equation*}
where$ \| \bm{\tau} \|^2_{D^{-1}} = (\bm{\tau}, \bm{\tau})_{D^{-1}} := (D^{-1}\bm{\tau}, \bm{\tau})$, $\forall \ \bm{\tau} \in \mathscr{W}$, and
$\| \bm{v} \|_L^2 = (\bm{v},\bm{v})_L := (L \bm{v},\bm{v})$, $\forall \ \bm{v} \in \mathscr{V}$.

In addition to the vertex space $\mathscr{V}$ and edge space $\mathscr{W}$, another important space of a graph~$\mathcal{G}$ is the so-called \emph{cycle space} (see \cite{bollobas2013modern} for more details), denoted by~$\mathscr{C}$, which is defined as,  
\begin{equation} \label{eqn: def_cycle_space}
\mathscr{C}:=\big\{ \bm{c} \in \mathscr{W} \ | \ G^\intercal \bm{c} = \bm{0} \big\}.
\end{equation}
Each cycle on the graph $\mathcal{G}$ corresponds to an element $\bm{c}$ in the cycle space $\mathscr{C}$. To be more specific, if $\{i_1,i_2,\ldots,i_k,i_{k+1}=i_1\}$ is a cycle $\mathcal{G}$, that is,
\[
  i_j\in\mathcal{V},\quad\mbox{and}\quad
  \left\{\max\{i_{j},i_{j+1}\},\min\{i_{j},i_{j+1}\}\right\}\in \mathcal{E},\quad j=1,\ldots,k,
\]
we define the components of $\bm{c}\in\mathscr{C}$ as follows:
\begin{equation*} \label{eqn: cycle_to_basis}
 ( \bm{c})_{e} =
  \operatorname{sign}(i_{j}-i_{j+1}), \quad
  \mbox{if}\quad e=
\left\{\max\{i_{j},i_{j+1}\},\min\{i_{j},i_{j+1}\}\right\},
\end{equation*}
with $(\bm{c})_e$ extended as zero for all edges that are not on the cycle. 
Besides its definition~\eqref{eqn: def_cycle_space}, we can also characterize the cycle space $\mathscr{C}$ by its basis. As discussed in~\cite{2009_CycleBasis}, the cycle space of a connected simple graph has dimension $m-n+1$  (by simple graphs we mean the graphs that contains neither self-loops nor duplicate edges). Note there is more than one way to find the basis of the cycle space (see the survey paper \cite{2009_CycleBasis}). For general graphs, a commonly used set of basis for the cycle space is the \emph{basis of fundamental cycles}. Such basis is not unique, but each such basis is induced by a spanning tree.  In a connected graph $\mathcal{G}$, with $n$ vertices, a spanning tree is any subgraph of $\mathcal{G}$ that connects all the vertices and has no cycles or, equivalently, has exactly $(n-1)$ edges. To construct the basis of fundamental cycles corresponding to a spanning tree $\mathcal{T} = (\mathcal{V}, \mathcal{E}_\mathcal{T}, \omega_\mathcal{T})$ of a graph~$\mathcal{G}$, we proceed as follows. For each edge that does not belong to the tree, i.e., $e  = \big\{ i,j \big\} \in \mathcal{E} \setminus \mathcal{E}_{\mathcal{T}}$, we can find a cycle $\big\{ i,j\big\} \cup p(i,j) $ where $p(i,j)$ is the path from vertex $i$ to vertex $j$ on the tree $\mathcal{T}$. Since the spanning tree $\mathcal{T}$ has exactly $(n-1)$ edges, there are $(m-n+1)$ such cycles. It can be shown that they are linearly independent~\cite{2009_CycleBasis} and, therefore, form a basis for the cycle space $\mathscr{C}$. 

\begin{figure}[h]
\centering
\begin{subfigure}[t]{.23\textwidth}
\includegraphics[width=\textwidth]{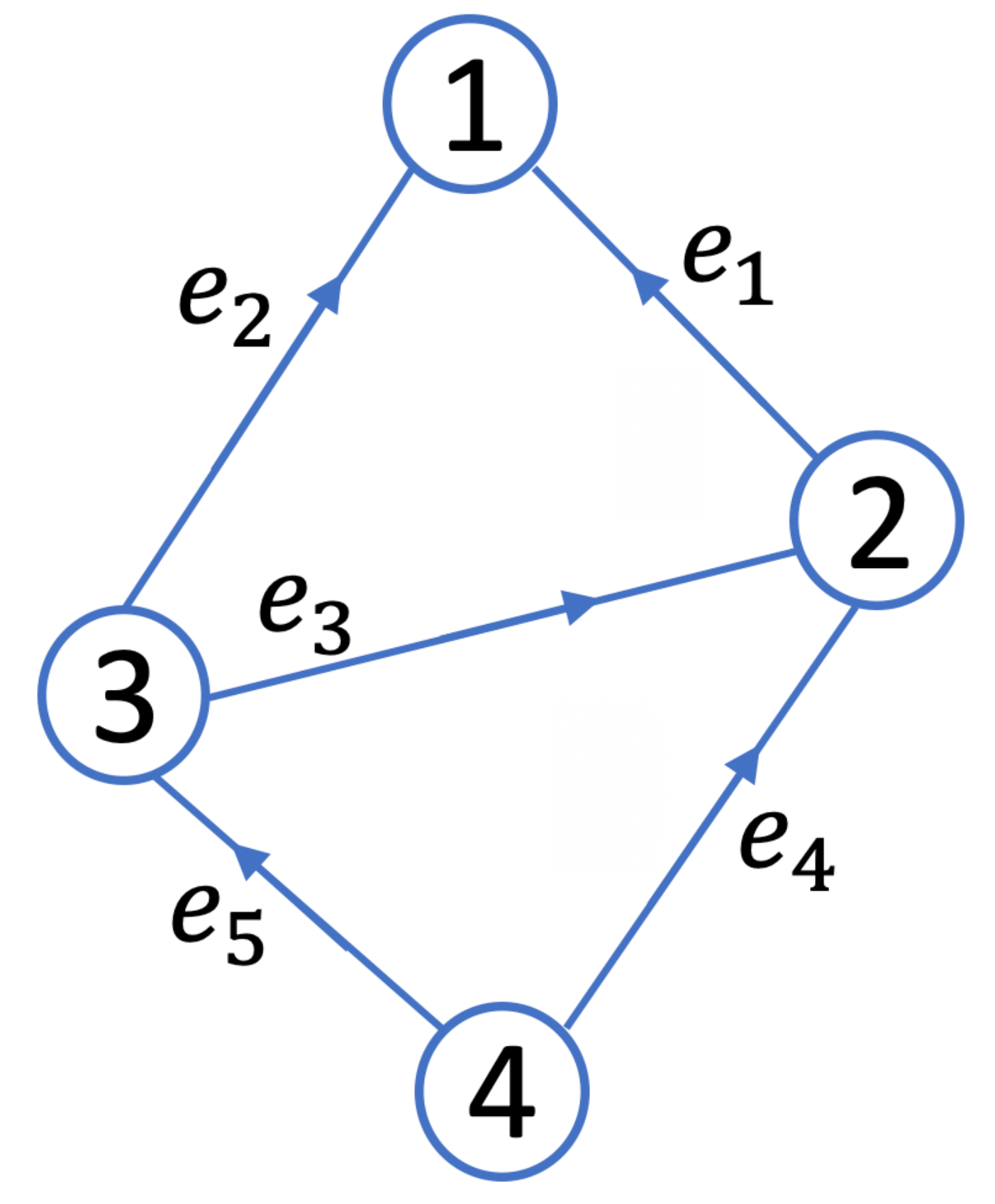}
\subcaption{original graph}
\end{subfigure}
\begin{subfigure}[t]{.23\textwidth}
\includegraphics[width=\textwidth]{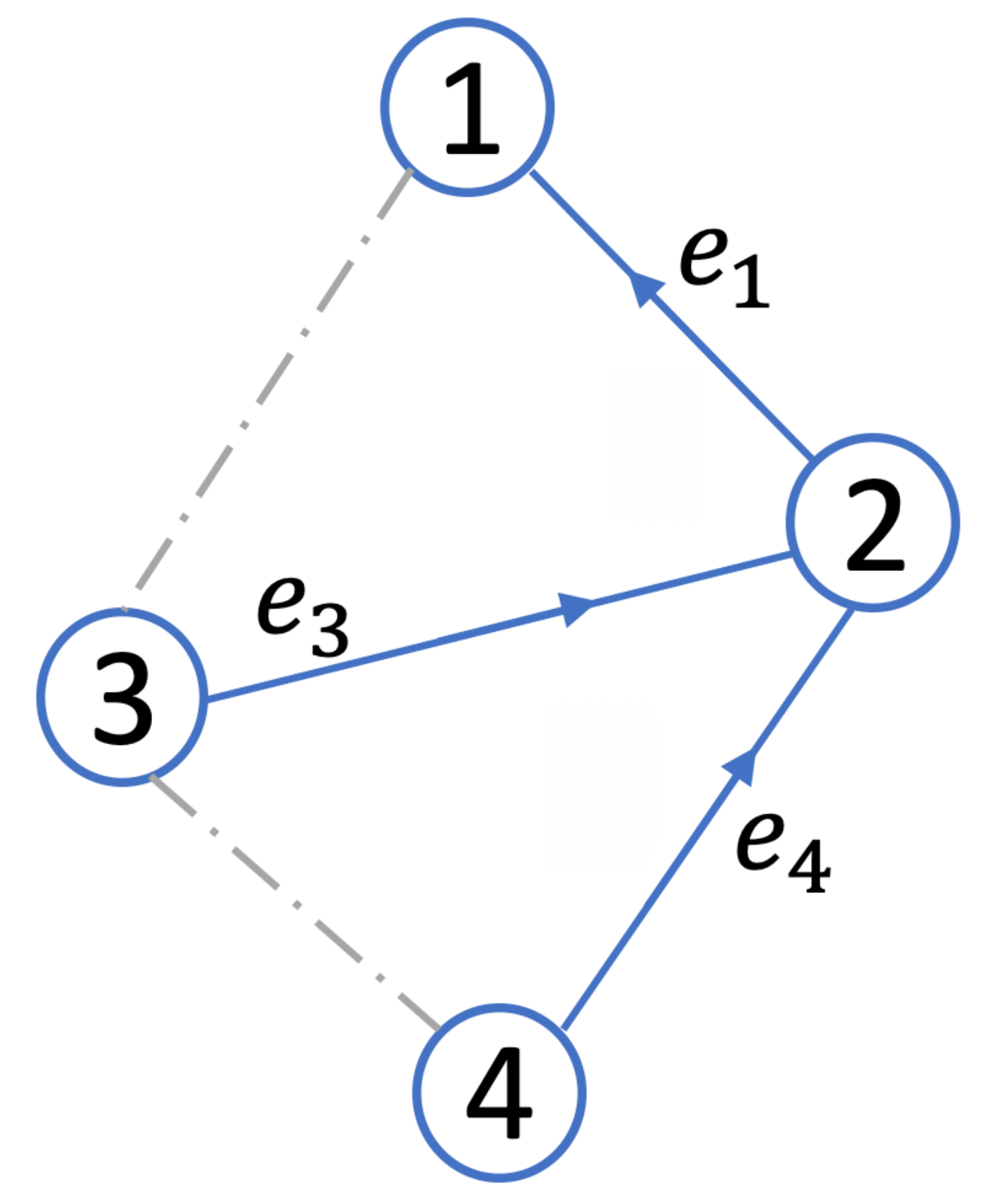}
\subcaption{spanning tree}
\end{subfigure}
\begin{subfigure}[t]{.23\textwidth}
\includegraphics[width=\textwidth]{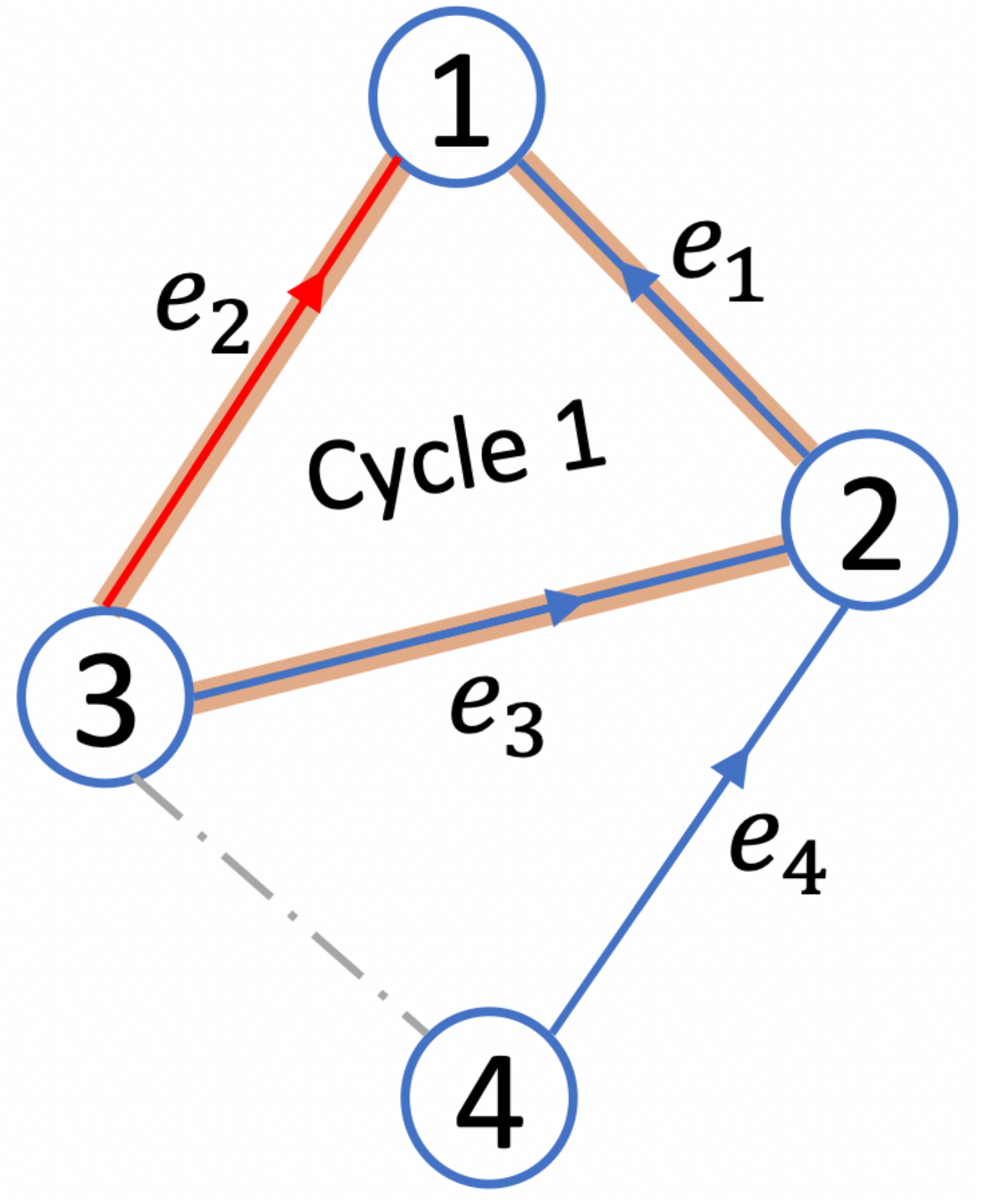}
\subcaption{adding edge $e_2$ to get cycle 1. }
\end{subfigure}
\begin{subfigure}[t]{.23\textwidth}
\includegraphics[width=\textwidth]{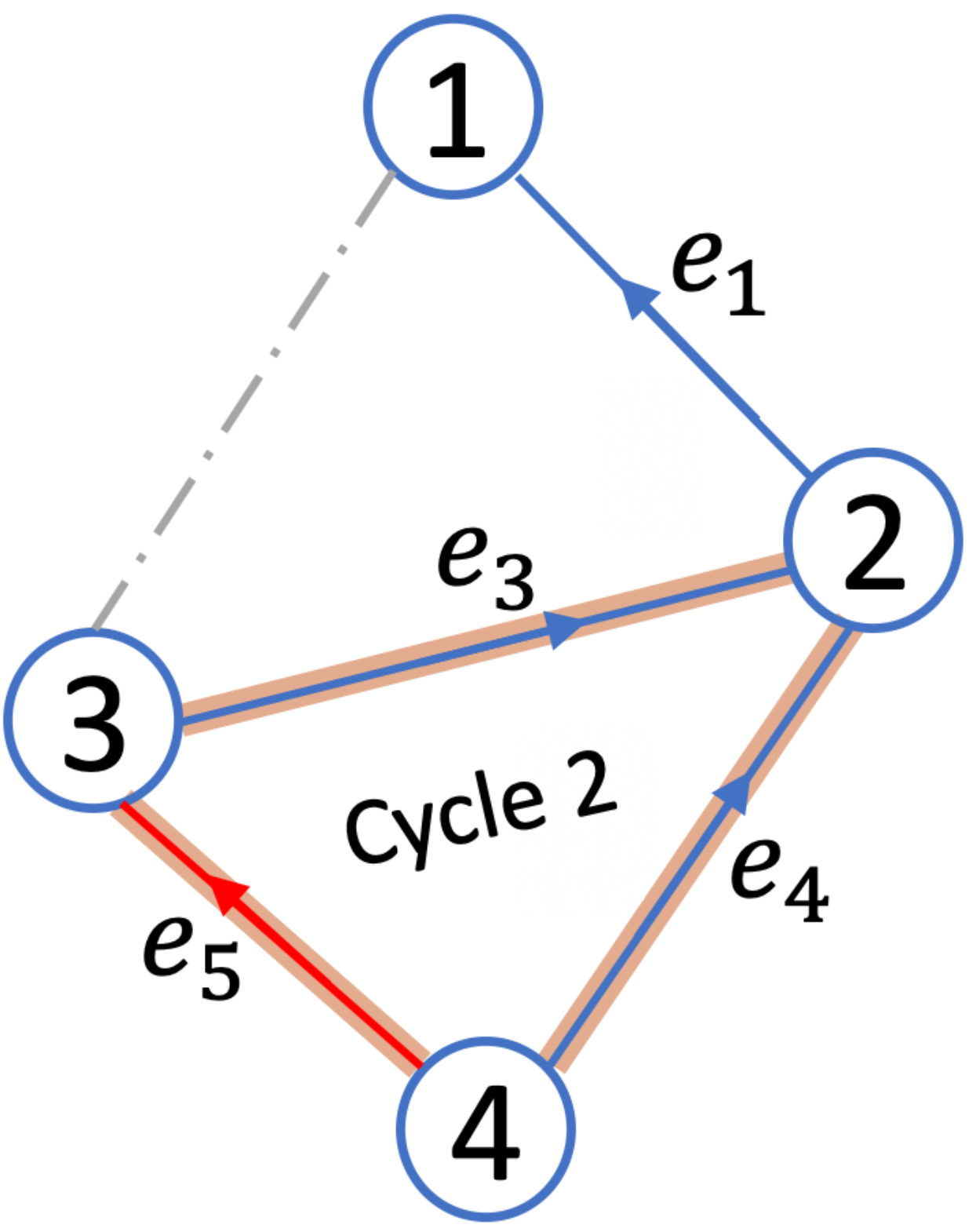}
\subcaption{add edge $e_5$ to get cycle 2.}
\end{subfigure}
\caption{Fundamental cycle basis}
\label{fig: cycle_basis}
\end{figure}

In \cref{fig: cycle_basis}, we give a simple example of the fundamental cycle basis. The tree in~\cref{fig: cycle_basis}(b) is a spanning tree of the graph in~\cref{fig: cycle_basis}(a). First, the edge $e_2$ is added back (see~\cref{fig: cycle_basis}(c)) which results in the first cycle $\left\{2,1,3,2 \right\}$ consisting of edges $e_1, e_2,$ and $e_3$.  The vector representation of the cycle induced by adding back edge $e_2$ is given by $ \bm{c}^{e_2} = [1,-1,1,0,0]^\intercal$.  Similarly, by adding edge $e_5$ back, we have the second cycle $\left\{ 2,3,4,2 \right\}$ formed by edges $e_3, e_4$ and $e_5$, which is represented by~$\bm{c}^{e_5} = [0,0,-1,-1,1]^\intercal$. $\bm{c}^{e_2}$ and $\bm{c}^{e_5}$ form a cycle basis. 

\subsection{Previous Results on A Posteriori Error Estimators}
We are interested in solving the following linear system of graph Laplacians:
\begin{equation} \label{eqn: lu=f}
L\bm{u}= \bm{f},
\end{equation}
by some iterative methods. Equation~\eqref{eqn: lu=f} simulates a rich spectrum of weighted graph Laplacians problems~\cite{Gutman72_EnergyNetworks, Hu2018_GraphLaplacian,Merris1994143,2012Gerta_ElectricalNetwork}. For example, in computational physics, the solution $\bm{u} \in \mathbb{R}^n$ models the potentials of an electrical flow where the resistance of each edge on the graph (electrical circuit) is the reciprocal of the edge weight, and $\bm{f} = [f_1, f_2, \cdots, f_n]^\intercal \in \mathbb{R}^n$ denotes the current supplied to each node $i$.  After $k$ iterations we get an approximated solution $\bm{u}^k$. If we can somehow construct the current error $ \bm{e}^k = \bm{u}-\bm{u}^k$, then the true solution will be easily obtained by $\bm{u} = \bm{u}^k+ \bm{e}^k$. In practice the true error $\bm{e}^k$ is not computable because $\bm{u}$ is unknown, so alternatively we seek to find $\widetilde{\bm{e}}^k$, an accurate estimation of $\bm{e}^k$, and use $\widetilde{\bm{e}}^k$ to improve the current approximation. Furthermore, an accurate estimation of the error gives an insight into the performance of the iterative methods.  For example, in AMG methods, such an estimation approximates the so-called ``smooth error'', which is responsible for the slow convergence of the AMG methods, and can be used to improve the AMG algorithm adaptively.  This leads to the adaptive AMG methods~\cite{aSA,2006Falgout_Adaptive, DAmbra.P;Vassilevski.P2015a,2012MacLachlanS_MoultonJ_ChartierT-aa,2008NagelA_FalgoutR_WittumG-aa} which has been actively researched over the past two decades.  

Since our a posteriori estimator is motivated by the a posteriori error estimator developed in~\cite{Xu2018_Agg}, we recall the main results and algorithms presented in~\cite{Xu2018_Agg} and start with the following fundamental lemma which relates the error and computed approximate solution. The proof of the lemma, as stated here, is found in~\cite{Xu2018_Agg}.
\begin{lemma} 
Let $\bm{u}$ be the solution to \eqref{eqn: lu=f}. Then for arbitrary $\bm{\tau} \in \mathscr{W}$, the following inequality holds for all $\bm{v}\in \mathscr{V}$:
\begin{equation} \label{eqn: prev_error_estimates}
\|\bm{u}-\bm{v}\|_L \leq \|DG\bm{v}- \bm{\tau}\|_{D^{-1}} + C_p^{-1} \|G^\intercal \bm{\tau}-\bm{f}\|_\mathscr{V}.
\end{equation}
where $C_p$ is Poincar\'e's constant of the graph Laplacian $L$.
\end{lemma}

For a fixed $\bm{v}$, denote the right-hand side of \eqref{eqn: prev_error_estimates} by:
 \begin{equation*}
 \eta(\bm{\tau}) =  \|DG\bm{v}-\bm{\tau}\|_{D^{-1}} + C_p^{-1} \|G^{\intercal} \bm{\tau}-\bm{f}\|_\mathscr{V}.
 \end{equation*}
This naturally provides a posteriori error estimator for estimating the error $\bm{u} -\bm{v}$ if $\bm{v}$ is an approximate solution, i.e., $\bm{v} = \bm{u}^k$.  Moreover, by minimizing the right-hand side of \eqref{eqn: prev_error_estimates} with respect to $\bm{\tau}$, we can obtain an accurate estimator. To solve the minimization problem efficiently, in~\cite{Xu2018_Agg}, an upper bound $E(\beta, \bm{\tau})$ of $\eta(\bm{\tau})$ was introduced as follows,
 \begin{equation*}
 \eta^2(\bm{\tau}) \leq E(\beta, \bm{\tau}),
 \end{equation*}
 where
 \begin{equation*}
 E(\beta, \bm{\tau}):= (1+\beta)\|DG\bm{v}-\bm{\tau}\|_{D^{-1}}^2 + (1+\frac{1}{\beta}) C_p^{-1} \|G^\intercal \bm{\tau}-\bm{f}\|_\mathscr{V}^2.
 \end{equation*}
 And an accurate estimator can be obtained by computing~$\min_{\beta, \bm{\tau}} E(\beta, \bm{\tau})$. In~\cite{Xu2018_Agg}, an alternating process is applied to minimize $E(\beta, \bm{\tau})$  with respect to $\beta$ (with the techniques proposed in~\cite{2011Tomar_MinBeta}) and $\bm{\tau}$ iteratively, as summarized in~\cref{alg:minimize_E} (see~\cite{Xu2018_Agg} for details):

\begin{algorithm}[htp!]
\caption{Alternating Process for Solving  $\min_{\beta, \bm{\tau}} E(\beta,\bm{\tau})$} \label{alg:minimize_E}
\begin{algorithmic}[1]
\Procedure{[$\beta, \bm{\tau}$] =MinimizeBound}{$\beta^0, \bm{\tau}^0$}
\For{$k = 1, 2, \cdots,$ max\_iter}
\State compute $\bm{\tau}^{k} = \mathrm{argmin}_{\bm{\tau}}E(\beta^{k-1}, \bm{\tau})$.
\State compute $\beta^{k} = \mathrm{argmin}_{\beta}E(\beta, \bm{\tau}^{k})$.
\EndFor
\EndProcedure
\end{algorithmic}
\end{algorithm}

Although the approach developed in~\cite{Xu2018_Agg} provides a reliable error estimator, the corresponding computational cost might be expensive due to the iterative minimization of $E(\beta, \bm{\tau})$ in step 3 and 4 in~\cref{alg:minimize_E}. One iteration of step 4 can be as difficult as solving the original system, which makes this approach expensive computationally. In order to improve the accuracy of the a posteriori error estimator, and, more importantly, to improve the efficiency of computing it, we develop a novel technique for estimating the error based on~\eqref{eqn: prev_error_estimates}, which we present next.

\section{Efficient Algorithm for Computing A Posteriori Error Estimator} \label{sec: main} 
As we have pointed out, the derivation of the a posteriori error estimator is based on approximating the  Helmholtz decomposition of the error. This provides a tighter error bound than the one proposed in~\cite{Xu2018_Agg} and can be implemented efficiently. 

\subsection{Hypercircle Identity and Error Estimation on Graphs}
Our design of an a posteriori error estimator is motivated by~\eqref{eqn: prev_error_estimates}.  For a given $\bm{f} \in \mathscr{V}$, we define the space $\mathscr{W}(\bm{f}) = \big\{ \bm{\tau} \in \mathscr{W} \;| \; G^\intercal \bm{\tau} =\bm{f} \big\}$.  If we choose $\bm{\tau} \in \mathscr{W}(\bm{f})$, then the second term on the right hand side of~\eqref{eqn: prev_error_estimates} vanishes and we only have the first term left. If we minimize this term with respect to $\bm{\tau} \in \mathscr{W}(\bm{f})$, we can immediately get an accurate estimation.  We summarize this in the following theorem.
\begin{theorem}  \label{thm:error}
Let $\bm{u}$ be the solution to \eqref{eqn: lu=f}. Then for any $\bm{v} \in \mathscr{V}$, we have,
\begin{equation} \label{eqn: error_estimates_exact}
\|\bm{u}-\bm{v}\|_L = \min_{\bm{\tau} \in \mathscr{W}(\bm{f})} \|DG\bm{v} - \bm{\tau}\|_{D^{-1}}.
\end{equation}
\end{theorem}

To prove~\cref{thm:error}, we will make use of the next lemma, also known as \emph{hypercyrcle identity} (see \cite{1947Prager_ErrorIdentity}):
\begin{lemma} \label{lem: error_identity}
Let $\bm{u}$ be the solution to \eqref{eqn: lu=f}. Then for any $\bm{v}\in \mathscr{V}$ and any $\bm{\tau}\in \mathscr{W}(f)$ the following identity holds:
\begin{equation*}
    \| \bm{u}-\bm{v}\|_L^2+\|DG\bm{u}-\bm{\tau}\|^2_{D^{-1}}=\|DG\bm{v}-\bm{\tau}\|^2_{D^{-1}}.
\end{equation*}
\end{lemma}
Now we are ready to prove~\cref{thm:error}.
\begin{proof}(Theorem 3.1)
 We first show $\|\bm{u}-\bm{v}\|_L \leq \min_{\bm{\tau} \in \mathscr{W}(\bm{f})} \|DG\bm{v} - \bm{\tau}\|_{D^{-1}}$. 
It follows from \cref{lem: error_identity} that,
\begin{equation*}
\|\bm{u}-\bm{v}\|_L^2 =  \|DG\bm{v}-\bm{\tau}\|_{D^{-1}}^2 - \|DG\bm{u}-\bm{\tau}\|_{D^{-1}}^2 \leq \|DG\bm{v}-\bm{\tau}\|_{D^{-1}}^2.
\end{equation*}
Since the inequality holds for any $\bm{\tau} \in \mathscr{W}(\bm{f})$, we have:
\begin{equation*}
\|\bm{u}-\bm{v}\|_L^2 \leq \min_{\bm{\tau} \in \mathscr{W}(\bm{f})} \|DG\bm{v} - \bm{\tau}\|_{D^{-1}}^2.
\end{equation*}
To show the other direction, note that, 
\begin{align*}
\|\bm{u}-\bm{v}\|_L^2 &= (L(\bm{u}-\bm{v}),\bm{u}-\bm{v}) = (G^\intercal DG (\bm{u}-\bm{v}), \bm{u}-\bm{v}) = (DG(\bm{u}-\bm{v}), G(\bm{u}-\bm{v})) \\
&=(DG(\bm{u}-\bm{v}),DG(\bm{u}-\bm{v}))_{D^{-1}} = \|DG(\bm{u}-\bm{v})\|_{D^{-1}}^2 = \|DG\bm{v}-DG\bm{u}\|_{D^{-1}}^2\\
&\geq \min_{\bm{\tau} \in \mathscr{W}(\bm{f})} \|DG\bm{v} - \bm{\tau}\|_{D^{-1}}^2.
\end{align*}
This completes the proof.
\end{proof}

Form~\cref{thm:error}, we observe that, 
\begin{equation*}
    \|\bm{u}-\bm{v} \|_{L} \leq \|DG\bm{v}-\bm{\tau} \|_{D^{-1}},
\end{equation*}
for any $\bm{\tau}\in\mathscr{W}(\bm{f})$. This motivates us to define the following computable quantity,
\begin{equation} \label{eqn: error_estimates}
\psi(\bm{\tau}) := \|DG\bm{v} - \bm{\tau}\|_{D^{-1}}, \quad \forall \, \bm{\tau} \in \mathscr{W}(f).
\end{equation}
If $ \bm{v} $ is the approximate solution to \eqref{eqn: lu=f}, $\psi(\bm{\tau})$ gives an a posteriori estimator for the true error $\bm{u}-\bm{v}$ for any choice of $\bm{\tau} \in \mathscr{W}(\bm{f})$. If $\bm{\tau}^*$ is the minimizer of the right-hand side of \eqref{eqn: error_estimates_exact}, then $ \psi (\bm{\tau^*}) = \|\bm{u}-\bm{v}\|_L$. Of course, computing the minimizer $\bm{\tau}^*$ exactly would be expensive. As an alternative we propose a Schwarz method for approximating $\bm{\tau}^*$. As our numerical tests show, we can obtain a reasonably good approximation $\bm{\tau} \in \mathscr{W}(\bm{f})$, $\bm{\tau}\approx\bm{\tau}^*$.

\subsection{Efficient Evaluation of the Error Estimator}

Our approach is to solve the minimization problem 
\begin{equation} \label{eqn: minimize_phi}
\min_{\bm{\tau} \in \mathscr{W}(\bm{f})} \|DG\bm{v} -\bm{ \tau}\|_{D^{-1}}
\end{equation}
 based on the Helmholtz decomposition of $\bm{\tau}$ on the graph, i.e.,   $\bm{\tau} = \bm{\tau}_f+\bm{\tau}_0^*$, where $\bm{\tau}_f \in \mathscr{W}(\bm{f})$ and $\bm{\tau}^*\in \mathscr{C}$ such that $(\bm{\tau}_f, \bm{\tau}_0^*) = 0$. Here $\bm{\tau}_f$ is curl-free and $\bm{\tau}_0^*$ is divergence-free. In particular, we first find a $\bm{\tau}_f \in \mathscr{W}(\bm{f})$ by solving a graph Laplacian on a spanning tree of the graph.  Then for a given $\bm{\tau}_f \in \mathscr{W}(\bm{f})$, the minimization problem \eqref{eqn: minimize_phi} becomes,
\begin{equation*}
\min_{\bm{\tau}_0 \in \mathscr{C}} \|DG\bm{v} - \bm{\tau}_f - \bm{\tau}_0\|_{D^{-1}}.
\end{equation*}
Solving this constrained minimization problem exactly will give the exact minimizer $\bm{\tau}_0^*$ and thus theoretically give the true error, which is an overkill in terms of finding a posteriori error estimator.  In practice, we only need to solve it approximately since as long as $\bm{\tau} \in \mathscr{W}(\bm{f})$, $\psi(\bm{\tau})$ will provide an upper bound of the error, which can be used as an error estimator. Note that this approximation is (always) subject to a trade-off: the error estimator will approximate the true error very accurately if we solve the optimization problem at unacceptable computation cost; or we can approximate and obtain a not-so-tight error estimator at an optimal computational cost. 

\subsubsection{Computing the Curl-free Component of the Error}
In this subsection, we discuss how to compute $\bm{\tau}_f \in \mathscr{W}(\bm{f})$ with optimal computational complexity for a given graph.  
For any $\bm{\tau}_f \in \mathscr{W}(\bm{f})$, we have,
\begin{equation} \label{eqn: compute_tauf}
G^\intercal \bm{\tau}_f  = \bm{f}.
\end{equation}
Since $G^\intercal$ is the discrete divergence operator on the graph, the solution to the above equation is not unique and difficult to compute in general.  However, we just need to find one $\bm{\tau}_f$.  Here, based on a spanning tree $\mathcal{T}$ of the graph $\mathcal{G}$, we present an approach with optimal complexity, i.e., $\mathcal{O}(n)$ computational cost. 

For a given spanning tree $\mathcal{T} = (\mathcal{V}, \mathcal{E}_\mathcal{T}, \omega_\mathcal{T})$, we look for a $\bm{\tau}_f$ satisfying~\eqref{eqn: compute_tauf} but only has nonzero entries on the edges that belong to the spanning tree $\mathcal{T}$.  In this case, we can rewrite~\eqref{eqn: compute_tauf} as:
\begin{equation} \label{eqn: break_down_tauf}
\bm{f} = G^\intercal \bm{\tau}_f = \begin{pmatrix}
G^\intercal_\mathcal{T} \; \; G^\intercal_{\mathcal{G} \setminus \mathcal{T}}
\end{pmatrix} \begin{pmatrix}
\bm{\tau}_{f, \mathcal{T}} \\ \bm{0}
\end{pmatrix},
\end{equation}
where $G^\intercal_\mathcal{T} $ is the discrete divergence operator that acts on edges in the tree $\mathcal{T}$, and $ G^\intercal_{\mathcal{G} \setminus \mathcal{T}}$ is the discrete divergence operator on edges that are only in the graph $\mathcal{G}$ but not in the tree $\mathcal{T}$.  From~\eqref{eqn: break_down_tauf},  we have, 
\begin{equation} \label{eqn: compute_tauf_tree}
G^\intercal_{\mathcal{T}} \bm{\tau}_{f, \mathcal{T}} = \bm{f}.
\end{equation}
Therefore, once we solve~\eqref{eqn: compute_tauf_tree}, we can assemble $\bm{\tau}_f$ by adding back the edges that are in the graph $\mathcal{G}$ but not in the tree $\mathcal{T}$.  Note that equation~\eqref{eqn: compute_tauf_tree} is defined only on the spanning tree $\mathcal{T}$.  We can first solve
\begin{equation} \label{eqn: lx=f_tree}
L_{\mathcal{T}} \bm{x} = \bm{f},
\end{equation}
$L_{\mathcal{T}}$ being the graph Laplacian of the tree $\mathcal{T}$. With the fact that 
 $L_{\mathcal{T}} = G^\intercal_{\mathcal{T}} D_{\mathcal{T}} G_{\mathcal{T}}$,
where $D_{\mathcal{T}}$ and $ G_{\mathcal{T}}$ are the discrete diagonal edge weight matrix and discrete gradient operator on tree $\mathcal{T}$, respectively, \eqref{eqn: lx=f_tree} can be rewritten as follows,
\begin{equation}\label{eqn: tree_laplacian}
G^\intercal_{\mathcal{T}} D_{\mathcal{T}} G_{\mathcal{T}} \bm{x} = \bm{f}.
\end{equation}
Comparing~\eqref{eqn: compute_tauf_tree} and~\eqref{eqn: tree_laplacian}, we naturally have,
\begin{equation*}
\bm{\tau}_{f, \mathcal{T}}:= D_{\mathcal{T}} G_{\mathcal{T}} \bm{x},
\end{equation*}
and can thereafter assemble the full $\bm{\tau}_f = (\bm{\tau}_{f, \mathcal{T}}, \bm{0})^{\intercal}$. The procedure for computing $\bm{\tau}_f$ is summarized in~\cref{alg:compute_tau_f}.  
\begin{algorithm}[htp!]
\caption{Computation of $\bm{\tau}_f$} \label{alg:compute_tau_f}
\begin{algorithmic}[1]
\Procedure{[$\bm{\tau}_f,\mathcal{T}$] =Compute$\bm{\tau}_f$}{$\mathcal{G}$, $\bm{f}$}
\State Build the spanning tree $\mathcal{T}$ from $\mathcal{G}$.
\State Solve $L_{\mathcal{T}}\bm{x} = \bm{f}$, where $L_{\mathcal{T}}=G^\intercal_{\mathcal{T}}D_{\mathcal{T}}G_{\mathcal{T}}$ .
\State Compute $\bm{\tau}_{f, \mathcal{T}} \gets D_{\mathcal{T}} G_{\mathcal{T}} \bm{x}$.
\State Assemble $\bm{\tau}_f$ as $\bm{\tau}_f \gets \begin{pmatrix} \bm{\tau}_{f, \mathcal{T}} \\ \bm{0} \end{pmatrix}  $.
\EndProcedure
\end{algorithmic}
\end{algorithm}

Identifying a spanning tree via the classic Breadth First Search algorithm (BFS tree) takes $\mathcal{O}(m+n)$ computational complexity for general graphs~\cite{cormen2009introduction}. If we consider sparse graphs in which $m = \mathcal{O}(n)$ or $\mathcal{O}(n \log \, n)$, this step costs at most $\mathcal{O}(n \log n)$.
The main computational cost of~\cref{alg:compute_tau_f} comes from Step 3, i.e., solving the linear system~\eqref{eqn: lx=f_tree}. As discussed in~\cite{1976RoseD_TarjanR_LuekerG-aa,2014VassilevskiP_ZikatanovL-aa}, it takes linear time to solve~\eqref{eqn: lx=f_tree}. 
Additionally, the matrix-vector multiplication in Step 4 has $\mathcal{O}(n)$ complexity for sparse graphs. Therefore, the overall complexity of~\cref{alg:compute_tau_f} is at most $\mathcal{O}(n\log n)$.

\subsubsection{Computing the the Div-free Component of the Error}
For a given $\bm{\tau}_f$, we need to solve the following constrained minimization problem,
\begin{equation} \label{eqn:opt-tau_0}
\bm{\tau}_0^* = \argmin_{\bm{\tau}_0 \in \mathscr{C}} \|DG\bm{v} - \bm{\tau}_f - \bm{\tau}_0\|_{D^{-1}}.
\end{equation}
The difficulty here is that we need to satisfy the constraint exactly when we compute an approximate $\bm{\tau}_0^*$ to get the estimated value of the error.  Our approach is to explicitly build the $(m-n+1)$ basis $\{ \bm{c}^e\}$ of the cycle space $\mathscr{C}$ as discussed in~\cref{sec:notations} and transform the constrained minimization problem~\eqref{eqn:opt-tau_0} into a unconstrained minimization problem. Denote by $\Omega$ the index set of the cycle basis. Then for any $\bm{\tau}_0 \in \mathscr{C}$ we can write $\bm{\tau}_0$ as a linear combination of the cycle basis: $\bm{\tau}_0 = \sum_{e\in \Omega}  \alpha_e \bm{c}^e$. Denote $\bm{\alpha} \in \mathbb{R}^{m-n+1}$, $(\bm{\alpha})_e = \alpha_e$. Thus, the minimization problem~\eqref{eqn:opt-tau_0} becomes,
\begin{equation} \label{eqn: minimization_problem}
\min_{\bm{\tau}_0 \in \mathscr{C}} \|DG\bm{v} - \bm{\tau}_f - \bm{\tau}_0\|_{D^{-1}}=\min_{\bm{\alpha}} \psi(\mathbb{\bm{\alpha}}) = \min_{\alpha_e, \ e\in \Omega} \| DG\bm{v}- \bm{\tau}_f - \sum_{e\in \Omega}  \alpha_e \bm{c}^e\|_{D^{-1}}.
\end{equation} 
This is an unconstrained least-squares problem and we can solve it with usual approaches. Moreover, the approximate solution is guaranteed to belong to the cycle space. Solving~\eqref{eqn: minimization_problem} exactly will eventually give us the exact error $\|\bm{u}-\bm{v}\|_L$. This step, however, has  a computational complexity comparable to  solving the original problem \eqref{eqn: lu=f}. Therefore, we solve it approximately via a few steps of one level overlapping Schwarz method (see~\cite{2002Xu_SubspaceCorrection, 2020Hu_SubspaceCorrection}). To describe such method, we first decompose the cycle space $\mathscr{C}$ into the following subspaces:
\begin{equation}
    \mathscr{C} = \mathscr{C}_1+ \mathscr{C}_2 + \cdots + \mathscr{C}_J.
\end{equation}
Note that $\mathscr{C}_i \cap \mathscr{C}_j$ is not necessarily empty.
Then, each step of the Schwarz method corresponds to a loop over all subspaces and solving the minimization problem \eqref{eqn: minimization_problem} in each of the subspace $\mathscr{C}_i$. That is, for $i = 1, 2, \cdots, J$, compute:
\begin{equation} \label{eqn: subdomain_minimization}
    \min_{\Delta \bm{\tau} \in \mathscr{C}_{i+1}} \|DG\bm{v} - \bm{\tau}_f - (\bm{\tau}^{i}_0+\Delta \bm{\tau} )\|_{D^{-1}}, 
\end{equation}
where $\bm{\tau}^i_0$ is the approximation to $\bm{\tau}_0^*$ after solving~\eqref{eqn: subdomain_minimization} in the first $i$ subspaces, and~$ \bm{\tau}_0^* = \mathrm{argmin}_{\bm{\tau}_0\in \mathscr{C}} \; \|DG\bm{v} - \bm{\tau}_f - \bm{\tau}_0 \|_{D^{-1}} $. The step of this relaxation method is completed after the approximation in $\mathscr{C}_J$ is computed. The solution obtained by iteratively solving equation~\eqref{eqn: subdomain_minimization} converges to the solution of the original minimization problem~\eqref{eqn:opt-tau_0} since this can be interpreted as a subspace optimization method whose convergence property was discussed in~\cite{xu1992iterative}.

 To keep low computational cost in computing the error estimator, we only run $\mathcal{O}(1)$ iterations of Schwarz method to approximately solve~\eqref{eqn: subdomain_minimization} for $\bm{\tau}_0^*$. Later in section~\ref{sec:num} we show that the error estimator computed with such an approximation is indeed accurate enough to capture the true error. The steps to compute $\bm{\tau}_0^* \in \mathscr{C}$ approximately are summarized in Algorithm~\ref{alg:compute_tau_0}, and we denote the approximation of $\bm{\tau}_0^*$ by $\bm{\tau}_0$.

\begin{algorithm}[htp!]
\caption{Computing an Approximation to $\bm{\tau}_0^*$} \label{alg:compute_tau_0}
\begin{algorithmic}[1]
\Procedure{$\bm{\tau}_0$ =Compute$\bm{\tau}_0$}{$\mathcal{G}, \mathcal{T}, \bm{\tau}_f$}
\State  Build the cycle basis $\big\{ \bm{c}^e\big\} $
\State Given initial guess $\bm{\tau}_0^{0}=\bm{0}$, 
\For{$i = 1,2,\cdots,$ max\_iter}  
\For{$k = 1, 2, \cdots, J$ } \Comment{iterate over each subdomain}
     \State $\Delta \bm{\tau}^* = \mathrm{argmin}_{\bm{\eta} \in \mathscr{C}_k} \;  \|DG\bm{v} - \bm{\tau}_f + \bm{\tau}_0^{(i-1)J+k-1} + \bm{\eta}\|_{D^{-1}} $. 
     \State $\bm{\tau}_0^{(i-1)J+k} = \bm{\tau}_0^{(i-1)J+k-1}+\Delta\bm{\tau}. $
\EndFor
\EndFor
\State \Return $\bm{\tau}_0^{J\cdot \text{max\_iter}}$.
\EndProcedure
\end{algorithmic}
\end{algorithm}

In~\cref{alg:compute_tau_0}, the cost of one step of the Schwarz method depends on the number of subspaces $J$ and the cost of solving~\eqref{eqn: subdomain_minimization} in each subspace. Here, we choose the following overlapping subspace decomposition: the $i$-the subspace is the span of the basis for the cycles incident with the vertex $i$. Thus, we have $J=n$, and
\begin{equation} \label{eqn:subspace_decomposition}
\begin{split}
  &\mathscr{C}_i = \operatorname{span}\{\bm{c}^j| \ \text{ cycle $j$ incidents with vertex $i$} \},  \quad i=1,\ldots,J.
\end{split}
\end{equation}
Since there are $n$ vertices we have $J = n$ subspaces. For sparse graphs using the special data structure proposed in~\cite{kelner2013simple}, the solution of the minimization problem~\eqref{eqn: subdomain_minimization} on each subspace $\mathscr{C}_i$ will be obtained with computational complexity at most~$\mathcal{O}(\log \, n)$. This holds even in the rare cases
when the dimension of $\mathscr{C}_i$ is large (e.g. $\sim n$). As a consequence the overall computational cost of each iteration of Schwarz method is $\mathcal{O}(n \log \, n)$ for this choice of subspace decomposition~\eqref{eqn:subspace_decomposition} and this assures a low computational cost in computing the proposed estimator.

\subsubsection{An Algorithm for A Posteriori Error Estimation with Helmholtz Decomposition}
Now we are ready to present the overall~\cref{alg:err_estimator} to (approximately) solve the minimization problem~\eqref{eqn: minimization_problem} and compute a posteriori error estimation for solving the graph Laplacian~\eqref{eqn: lu=f}.

\begin{algorithm}[H]
\caption{Computtion of the error estimator  $\min_{\bm{\tau}_0 \in \mathscr{C}} \|DG\bm{v} - \bm{\tau}_f - \bm{\tau}_0\|_{D^{-1}}$ }\label{alg:err_estimator}
\begin{algorithmic}[1]
\Procedure{$\psi$ =ErrorEstimates}{$\mathcal{G},\bm{v}$, $\bm{f} $}
\State [$\bm{\tau}_f, \mathcal{T}$] = Compute$\bm{\tau}_f$($\mathcal{G}, \bm{f}$).
\State $\bm{\tau}_0$ = Compute$\bm{\tau}_0$($\mathcal{G}, \mathcal{T}, \bm{\tau}_f$).
\State$\psi\gets \|DG\bm{v} - \bm{\tau}_f - \bm{ \tau}_0\|_{D^{-1}}$.
\Comment{Compute the value of the estimator}
\State \Return $\psi$.
\EndProcedure
\end{algorithmic}
\end{algorithm}

In~\cref{alg:err_estimator}, Step 2 to compute $\bm{\tau}_f$ takes $\mathcal{O}(n \log \, n)$ for any sparse graphs. Step 3 to compute $\bm{\tau}_0$ has complexity $\mathcal{O}(n \log \, n)$ for sparse graphs since the minimization problem~\eqref{eqn: minimization_problem} is solved approximately with $\mathcal{O}(1)$ iterations of Schwarz method. As a result, the overall computational complexity of~\cref{alg:err_estimator} is $\mathcal{O}(n \log \, n)$ for sparse graph $\mathcal{G}$. 

To make the a posteriori error estimator more useful, especially for developing the adaptive AMG methods for solving graph Laplacians~\cite{Xu2018_Agg, aAMG_graph1, Livne.O;Brandt.A.2012a}, we need to localize the a posteriori error estimator. Since,
\begin{equation*}
  \begin{aligned}
    \psi^2(\bm{\tau}) &= \|DG\bm{v} - \bm{\tau}\|_{D^{-1}}^2 = (DG\bm{v} - \bm{\tau})^T D^{-1}(DG\bm{v} - \bm{\tau})\\
    &= \sum_{e\in \mathcal{E}} \frac{1}{\omega_e}\big((DG\bm{v} - \bm{\tau})e \big)^2,
  \end{aligned}
\end{equation*}
we can localize the error estimator on each edge $e$ as follows,
\begin{equation} \label{eqn: local_error_estimator}
\psi_e({\bm{\tau}}) = \bm{\omega}_e^{-\frac{1}{2}} | (DG\bm{v}-\bm{\tau})_e |.
\end{equation}
We comment that the above localized error estimators is obtained for free in practice, since we have $(DG\bm{v}-\bm{\tau})_e$ available from the computation of the global error estimator $\psi(\bm{\tau})$ (see Step 4 in~\cref{alg:err_estimator}). 

This localized error estimator~\eqref{eqn: local_error_estimator} then can be used to design adaptive AMG methods. For example, it can be utilized in generating coarser aggregations that approximate the fine aggregates (vertices) accurately~\cite{Xu2018_Agg} or generate approximations to the level sets of the error for the path cover adaptive AMG method~\cite{Hu2018_PathCover}. 

\section{Numerical Results}\label{sec:num}

In this section, we present results of some numerical experiments demonstrating the efficiency of the a posteriori error estimator.

\subsection{Tests on 2D Uniform Grids}
We first test the performance of the algorithm on the unweighted graph Laplacian $L$ of 2D uniform triangle grids, which corresponds to solving a Poisson equation on a 2D square domain with Neumann boundary condition. The uniform triangle grid with grid size $h=2^{-l}$, $l = 5,6,7,8,9$ is used, and we take $\bm{u} = \sin(\frac{\pi}{2}\bm{x})\sin(\frac{\pi}{2}\bm{y})$. We set the approximate solution $ \bm{v}=\bm{0}$ and obtain the a posteriori error estimator $ \psi(\bm{\tau})$ with~\cref{alg:err_estimator}, in which the minimization problem \eqref{eqn: minimization_problem} is solved approximately with several iterations of the overlapping Schwarz method. We use the face cycle bases that correspond to the small triangles in the grid (cycle length is 3). With this choice of cycle basis, each of the decomposed subspaces in~\eqref{eqn:subspace_decomposition} have dimension $\mathcal{O}(1)$ since there are at most six cycles incident from a given vertex $i$. The low dimension of the subspaces assures that solving~\cref{eqn: subdomain_minimization} costs no more than $\mathcal{O}(1)$ computation and thus the computation cost of one iteration of Schwarz method remains $\mathcal{O}(n)$. 

In~\cref{Tab: global_bound}, we report the true error and the a posteriori error estimator~$\psi(\bm{\tau})$ on graph Laplacian systems of different scales.  $\displaystyle e_{\mathit{ff}}:=\frac{\psi (\bm{\tau})}{\|\bm{u}-\bm{v}\|_L} $ is also reported to show the efficiency of the error estimator. From~\cref{Tab: global_bound}, we observe that the CPU time for one iteration of Schwarz method grows linearly as the size of graph Laplacian systems increases. The error estimator $\psi(\bm{\tau})$ gradually approaches the true error $\|\bm{u}-\bm{v}\|_L$ when we increase the steps of Schwarz iteration.

\setlength{\tabcolsep}{5pt}
\begin{table}[H]
\centering
\caption{Efficiency of the error estimator on graph Laplacian systems on uniform triangle grids of different sizes. The value of the estimator $ \psi(\bm{\tau})$ is computed by solving \eqref{eqn: minimization_problem}  approximately with $1$, $3$, and $5$ iterations of the overlapping Schwarz method. The CPU time (in seconds) is also shown in the table.}
\label{Tab: global_bound}
\begin{tabular}{|r|c|c|c|c|c|c|c|c|c|c|}
\hline
      &            & \multicolumn{3}{c|}{1 iter}   & \multicolumn{3}{c|}{3 iters} & \multicolumn{3}{c|}{5 iters} \\ \hline
 \multicolumn{1}{|c|}{$|\mathcal{V}|$}   & $\|\bm{u}-\bm{v}\|_L$ & $ \psi(\bm{\tau})$  & $e_{\mathit{ff}}$    & time     & $ \psi(\bm{\tau})$  & $e_{\mathit{ff}}$      & time   & $ \psi(\bm{\tau})$  & $e_{\mathit{ff}}$      & time   \\ \hline
$1089$  & 1.73     & 2.25    & 1.30 & 0.03   & 1.99     & 1.15 & 0.04 & 1.91     & 1.10 & 0.06  \\ \hline
$4225$  & 1.73     & 2.67    & 1.55 & 0.05   & 2.28     & 1.32 & 0.11   & 2.16     & 1.25 & 0.16 \\ \hline
$16641$ & 1.73     & 3.36    & 1.95 & 0.14   & 2.76      & 1.60 & 0.37 & 2.56     & 1.48 & 0.62 \\ \hline
$66049$ & 1.72      & 4.43     & 2.57 & 0.53 & 3.51     & 2.03 & 1.40 & 3.20        & 1.86 & 2.31 \\ \hline
$263169$ & 1.72     & 6.01    & 3.49 & 1.92& 4.66     & 2.71 & 5.64 & 4.19    & 2.43 & 9.53 \\ \hline
\end{tabular}
\end{table}

More importantly, we would like to know whether the localized error estimator \eqref{eqn: local_error_estimator} approximates the true error on each edge accurately, since the localized estimation is the key to an effective coarsening scheme in adaptive AMG. Take $L$ as the weighted graph Laplacian of the uniform grid with grid size $h=2^{-5}$, $\bm{u} = \sin(\frac{\pi}{2}\bm{x})\sin(\frac{\pi}{2}\bm{y})$ and $ \bm{v}$ obtained by three iterations of Gauss Seidel method with random initial guess. We compute the error estimator using three iterations of the Schwarz method to solve the minimization problem in~\cref{alg:err_estimator}. 

\begin{figure}[ht!]
    \centering
    \includegraphics[width = 2.8in]{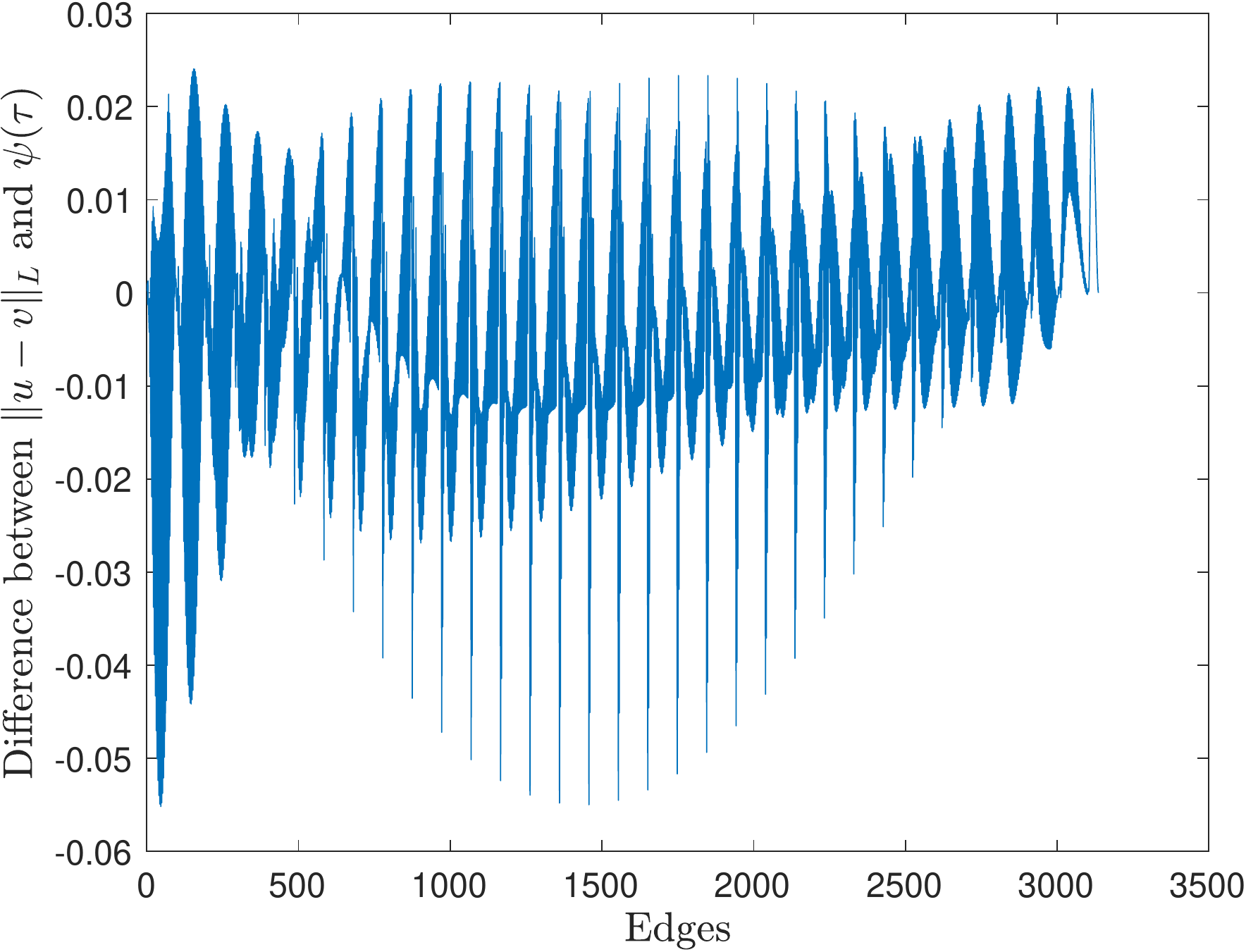}
    \caption{Difference between the true error $\|\bm{u}-\bm{v}\|_L $ and error estimator $\|DG\bm{v} - \bm{\tau} \|_{D^{-1}}$ on each edge $e$.}
    \label{fig: error_difference}
\end{figure}

In~\cref{fig: error_difference}, we plot the difference between the true error and the error estimator on each edge. On most of the edges the error estimator captures the true error well since the difference $\psi(\bm{\tau}) -\|\bm{u}-\bm{v}\|_L$ is no larger than $0.02$.

\subsection{Tests on ``Real World" Graphs}
In this section we test the proposed error estimator on some real world graphs from the SuiteSparse Matrix Collection \cite{RealWorldGraph}. We pre-process the undirected graphs by extracting the largest connected component of each graph and deleting self-loops. For each of these graphs, if the original edge weight is negative, we take its absolute value.

\begin{table}[htbp!]
\centering
\caption{Efficiency of error estimator on graph Laplacian systems arising from real world applications. The value of the estimator $\psi(\bm{\tau})$ is computed by solving \eqref{eqn: minimization_problem}  approximately with 3 iterations of Schwarz method. The graph types tested are unweighted (u) and weighted (w).}
\label{Tab: RealWorld}
\begin{tabular}{|r|r|r|l|c|c|c|c|}
\hline
\multicolumn{1}{|c|}{ID}   & 
\multicolumn{1}{c|}{$|\mathcal{V}|$} & 
\multicolumn{1}{c|}{$|\mathcal{E}|$} & 
\multicolumn{1}{c|}{Problem Type}   & 
\multicolumn{1}{c|}{Type} & 
\multicolumn{1}{c|}{$\|\bm{u}-\bm{v}\|_L$}  & 
\multicolumn{1}{c|}{$ \psi(\bm{\tau})$} & 
\multicolumn{1}{c|}{$e_{\mathit{ff}}$}  \\ \hline
8    & 292      & 958      & Least Squares Problem      & u        & 1.74     & 1.75  & 1.00  \\ \hline
1196 & 1879     & 5525     & Circuit Simulation & w        & 2.71     & 2.71  & 1.00 \\ \hline
22   & 5300     & 8271     & Power Network       & u        & 5.82     & 5.82  & 1.00 \\ \hline
1614 & 2048     & 4034     & Electromagnetic Problem   & w        & 0.47   & 0.50  & 1.07 \\ \hline
33   & 1423     & 16342    & Structural Problem         & w        & 14.5    & 19.7  & 1.36 \\ \hline
791  & 8205     & 58681    & Acoustic Problem          & w        & 23.8    & 37.7  & 1.58 \\ \hline
2777 & 1857& 13762 & Social Network & u & 52.9 & 76.3 & 1.44 \\ \hline
1533 & 2361 & 13828 & Protein Network& u& 4.61 & 4.70 & 1.01 \\ \hline 
\end{tabular}
\end{table}
\setlength{\tabcolsep}{\tabcolsepold}

\begin{figure}
    \centering
\begin{subfigure}[t]{0.32\textwidth}
\includegraphics[width = \textwidth]{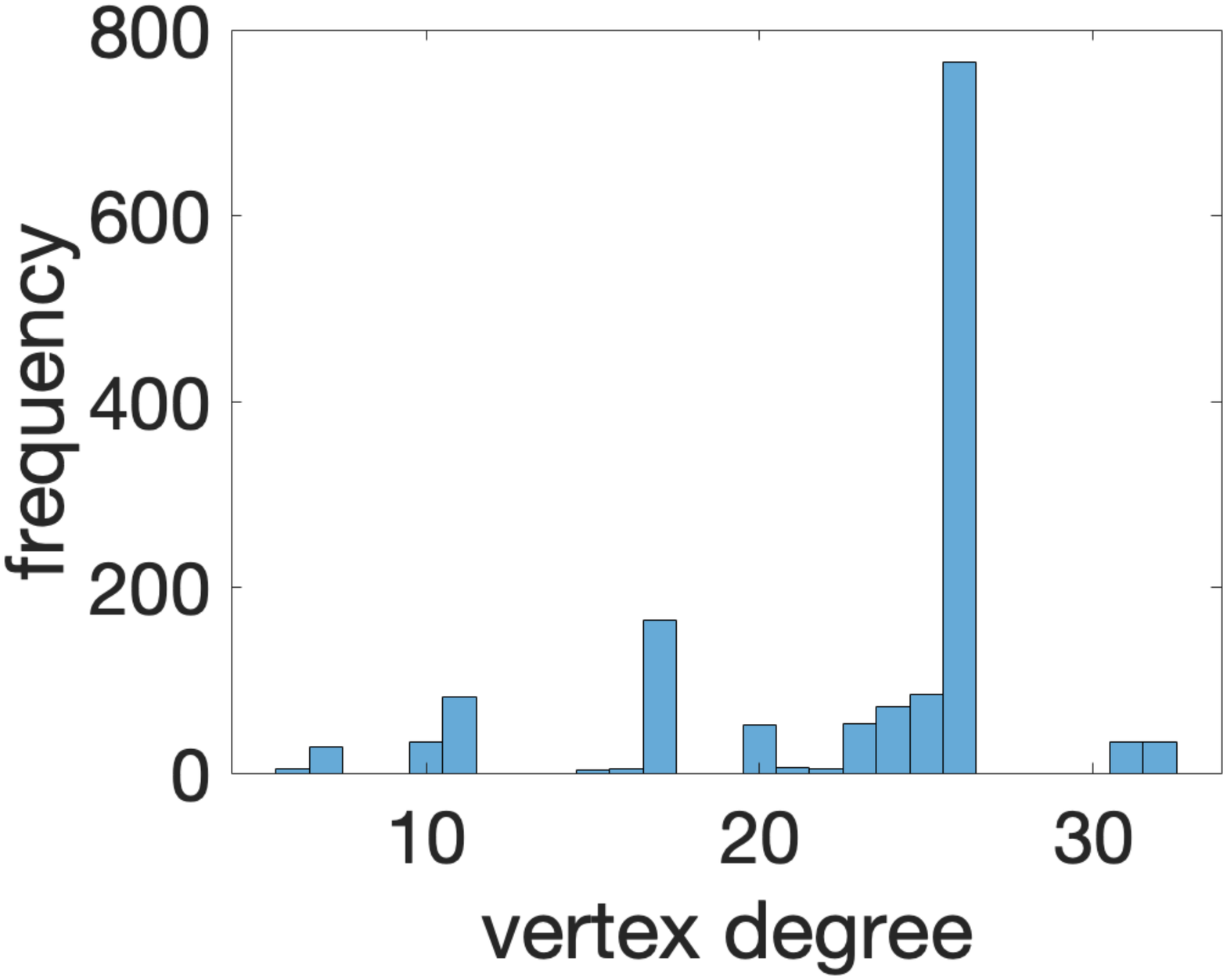}
\end{subfigure}
\begin{subfigure}[t]{0.32\textwidth}
\includegraphics[width = \textwidth]{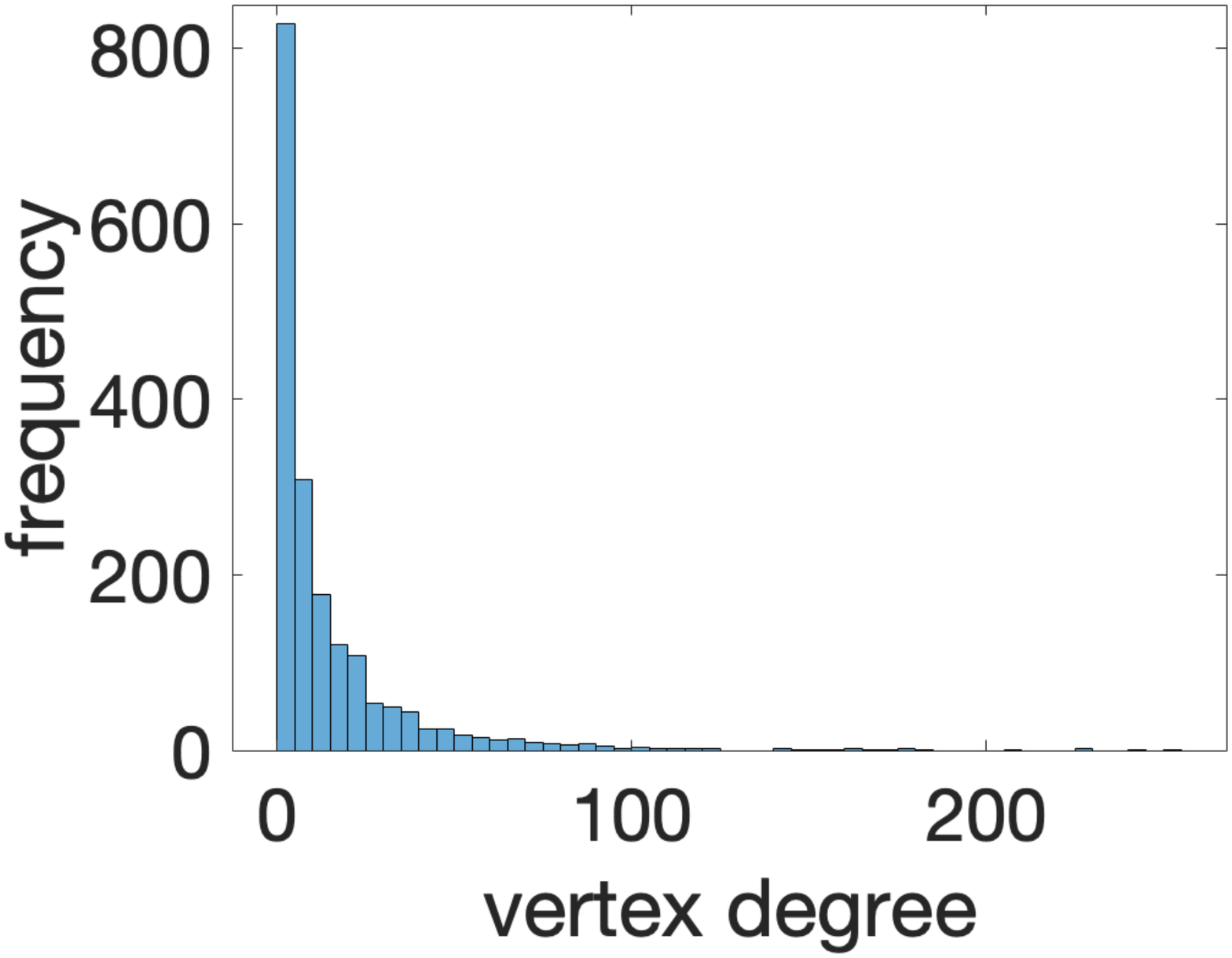}
\end{subfigure}
\begin{subfigure}[t]{0.32\textwidth}
\includegraphics[width = \textwidth]{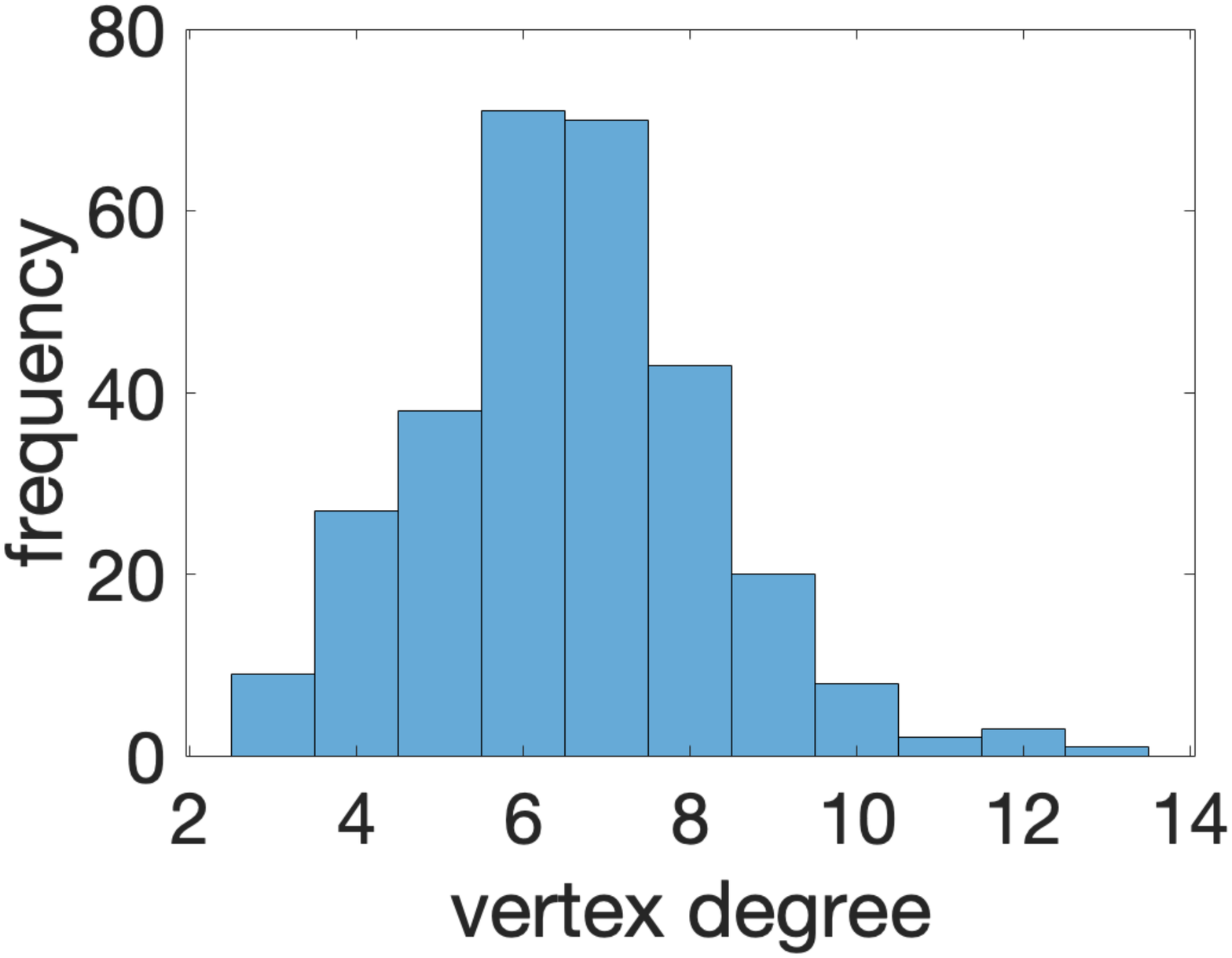}
\end{subfigure}
\caption{Representative degree distribution of networks studied in Table 2. \emph{left}:  Network ID 33. \emph{middle}:  Network ID: 2777 (power law distribution). \emph{right:} Network ID 8 (normal distribution).}
\label{fig: degree_distribution}
\end{figure}

In~\cref{Tab: RealWorld}, we summarize the basic information of the graphs and the performance of the error estimator. In our setting $\bm{u}$ is the exact solution for a problem with arbitrarily chosen right-hand side $\bm{f}$. The approximate solution $\bm{v}$ is obtained as a result of three iterations of the Gauss-Seidel method with this right-hand side. To compute the error estimator, we first use the breadth-first-search algorithm to find a spanning tree and then construct the spanning-tree-induced fundamental cycle basis. Then we apply three steps of Schwarz iterations to solve the minimization problem in~\cref{alg:err_estimator} to compute the overall error estimator. In \cref{fig: degree_distribution} we plot the degree distribution of selected networks. The differences in these degree distributions suggest that these networks have distinctive structural and dynamical properties.  As we can see from the results, for real-world graphs with different sizes, structures, and density, the error estimators approximate the true errors well in all cases, which demonstrate the effectiveness of our proposed algorithm for computing the a posterior error estimator. 

\subsection{Tests on Building Aggregations for AMG}\label{sec:PC-AMG}
Designing an effective coarsening scheme that projects the smooth error onto coarse levels accurately is the key to the AMG methods.  One example is the path cover adaptive algebraic multigrid (PC-$\alpha$AMG) proposed in~\cite{Hu2018_PathCover} for solving weighted graph Laplacian linear systems. Provided an accurate estimation of the current error, PC-$\alpha$AMG first forms vertex-disjoint path cover where paths approximate the level sets of the smooth error, then aggregates vertices to form aggregations. 

While an accurate error estimate makes PC-$\alpha$AMG highly efficient, the step to compute estimation of the true error in the original PC-$\alpha$AMG appears to be expensive in~\cite{Hu2018_PathCover}. Our a posteriori error estimator offers an ideal replacement which gives an accurate approximation in an efficient manner.

Here, we present the resulting path-cover aggregation obtained by using the proposed a posteriori error estimator for solving the unweighted graph Laplacian on the 2D uniform grid. We solve the corresponding linear system approximately with several iterations of relaxation scheme and curated the exact true error of the current solution (plotted in Figure 4(a)). Note that in practice this exact error is not directly accessible, and we approximate it with the proposed a posteriori error estimator, plotted in Figure 4(d). In Figure 4(b) and 4(e) we compare the path-cover aggregation generated with the true error and the a posteriori error estimator (computed approximately with 3 iterations of Schwarz method). The aggregation patterns are very similar.

To check whether the smooth error is transferred and represented accurately on the coarse level, we restrict the smooth error to the coarse level and then prolongate it back. We plot the differences between the smooth error and error after restriction and prolongation in Figure 4(c) for the case where the coarse level is constructed based on the exact smooth error, and in Figure 4(f) for the case where the coarse level is constructed based on a posteriori error estimator, respectively.  As we can see, although the shapes of differences in the two cases are different, the magnitudes of both cases are $0.045$, which indicates that the aggregation built with the error estimator is effective in capturing the true smooth error.  
 
\begin{figure}[H]
\centering
\begin{subfigure}[t]{.32\textwidth}
\includegraphics[width=\textwidth]{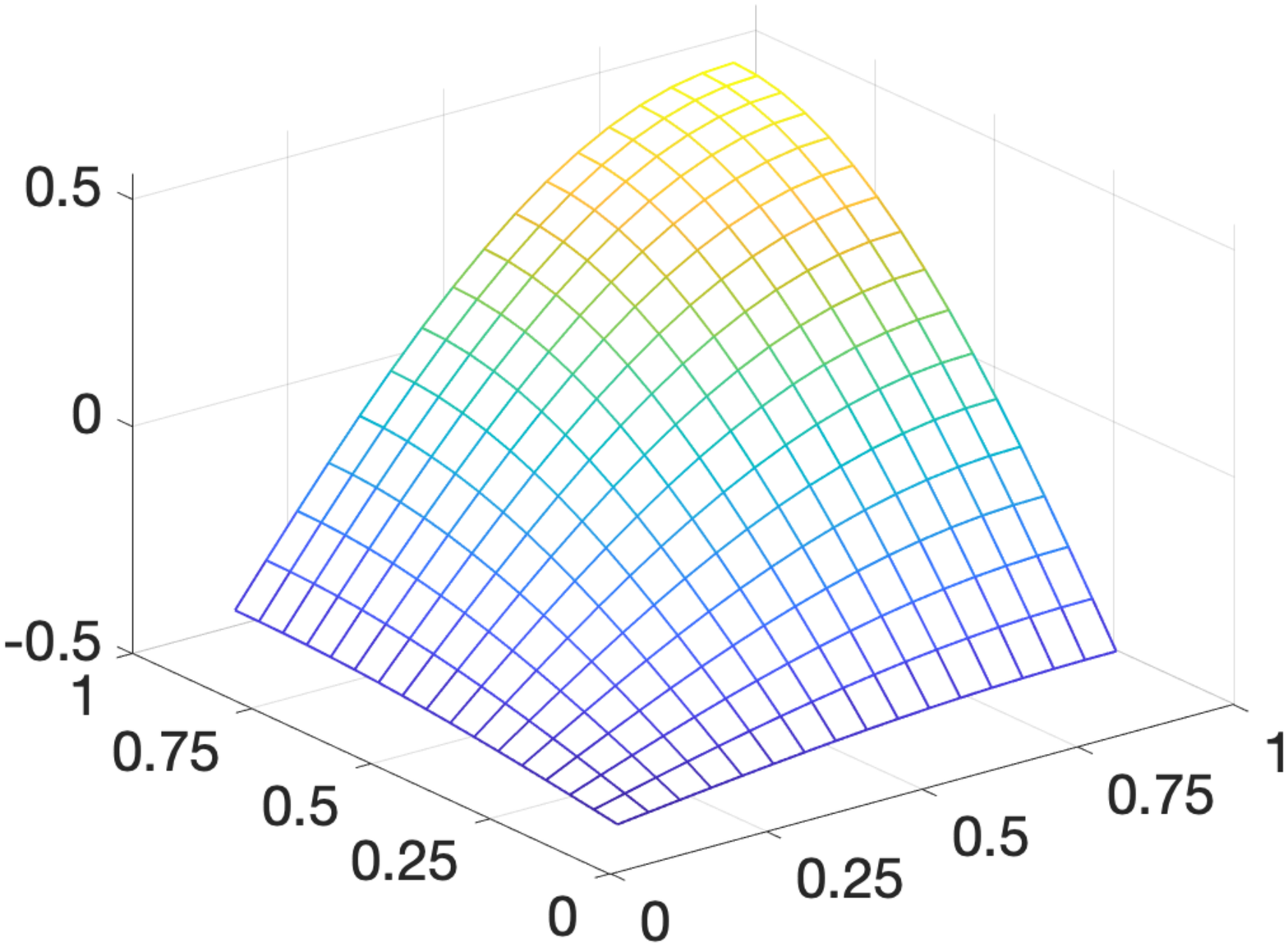}
\subcaption*{(a) True error of the current solution.}
\end{subfigure}
\begin{subfigure}[t]{.28\textwidth}
\includegraphics[width=\textwidth]{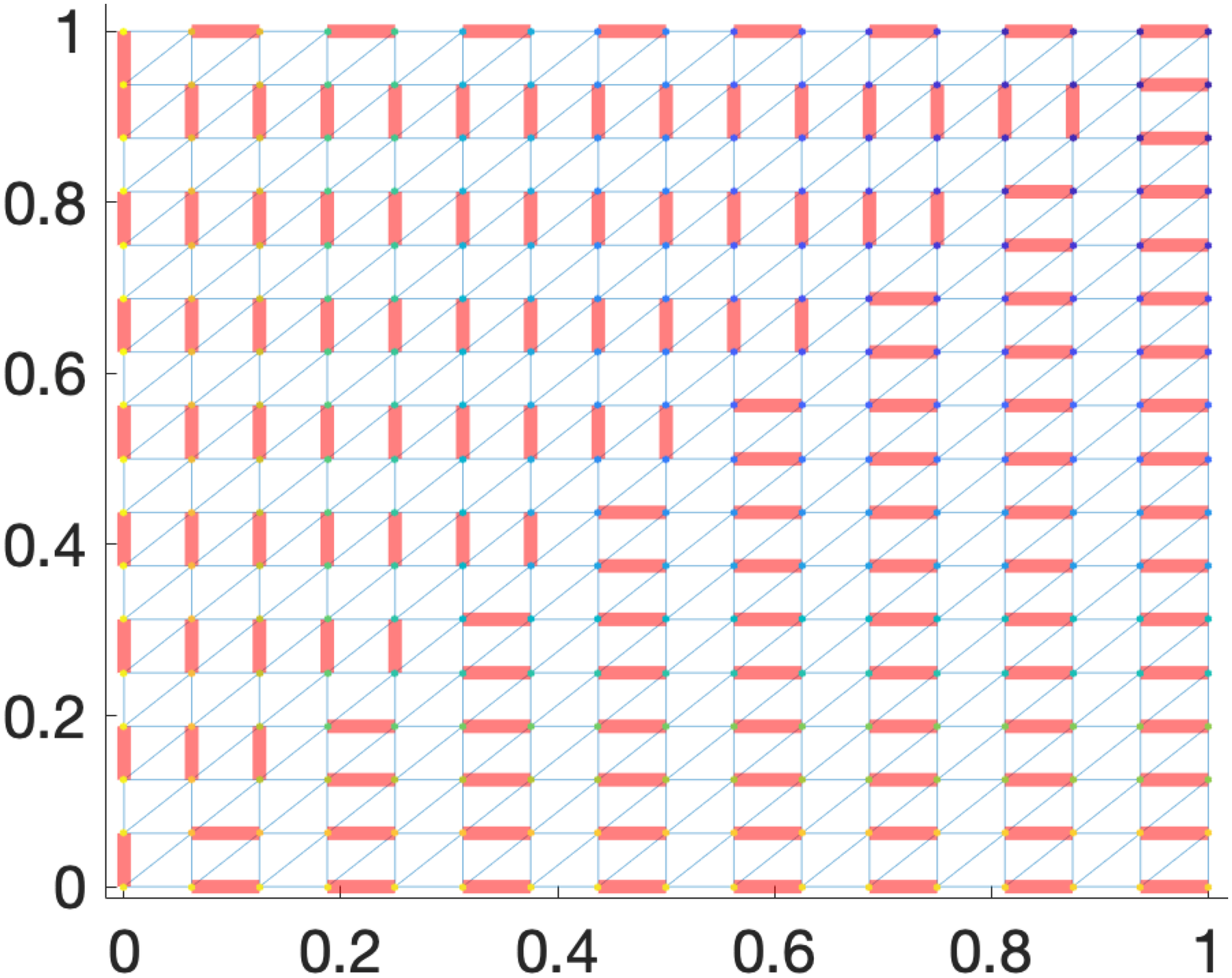}
\subcaption*{(b) Path cover aggregation based on the known exact error.}
\end{subfigure}
\begin{subfigure}[t]{.32\textwidth} 
\includegraphics[width=\textwidth]{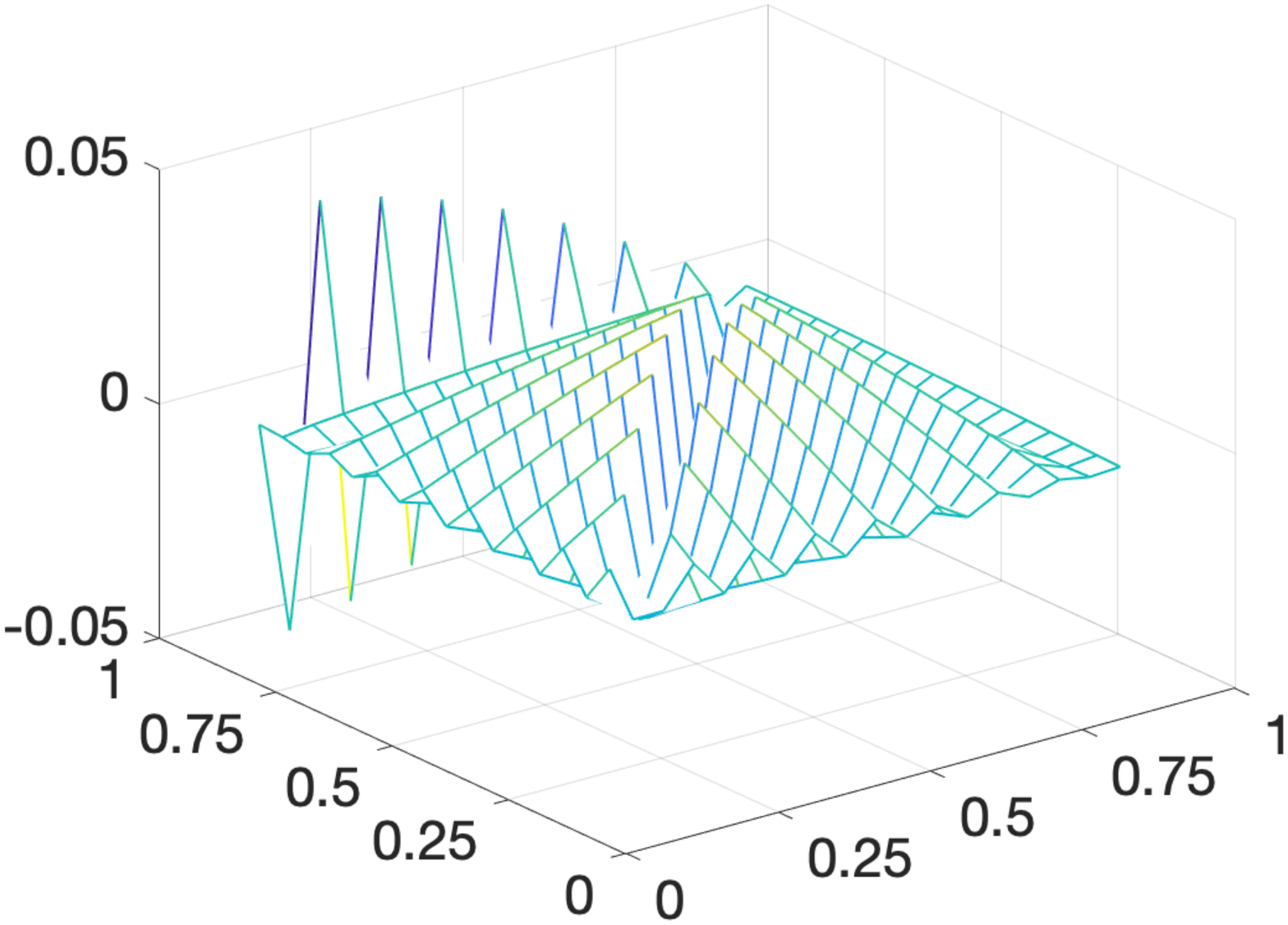}
\subcaption*{(c) Difference between the smooth error and the error after restriction and prolongation with aggregation based on true error.}
\end{subfigure}
\begin{subfigure}[t]{.32\textwidth}
\includegraphics[width=\textwidth]{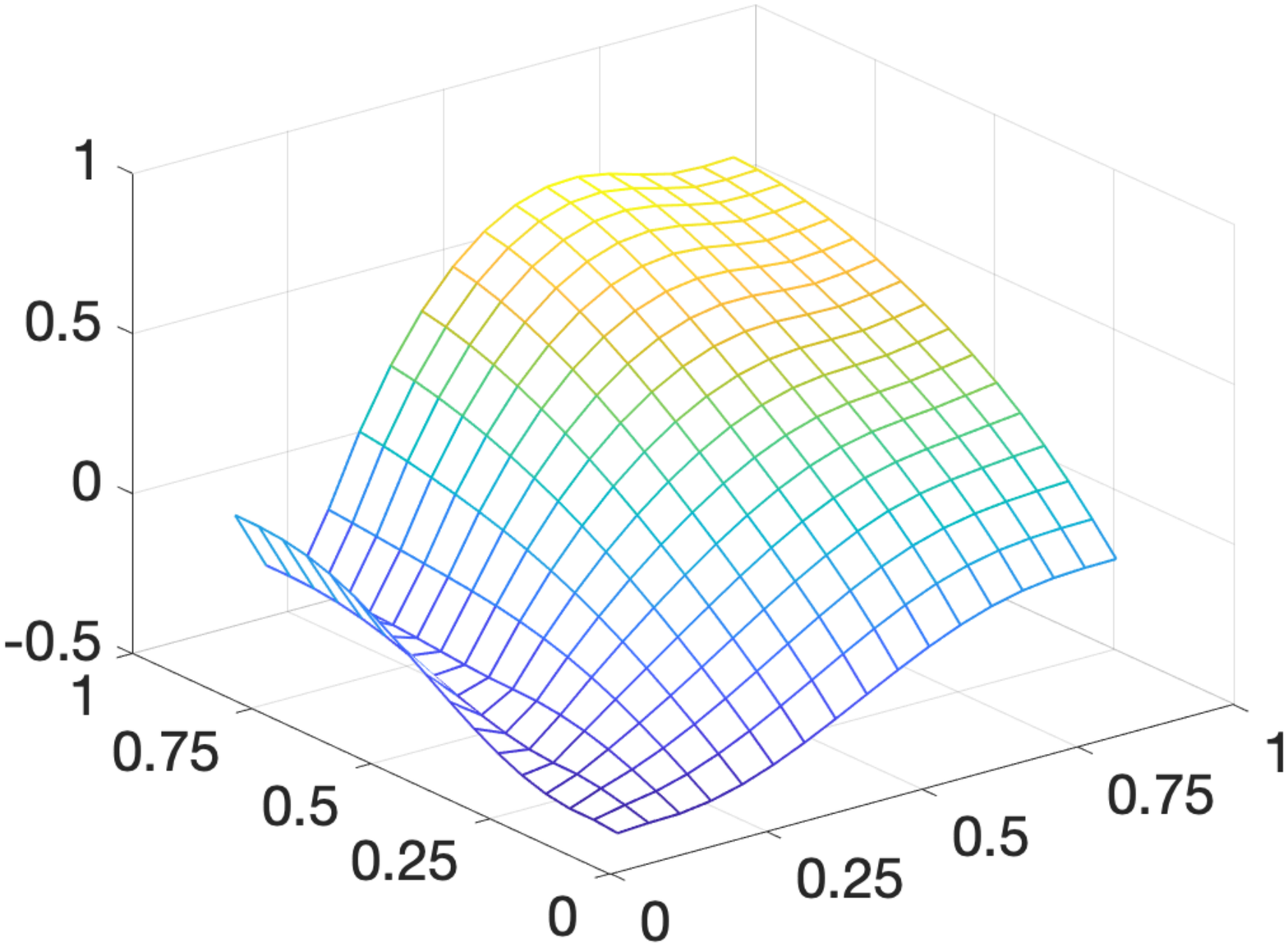}
\subcaption*{A posteriori error estimates of the current solution.}
\end{subfigure}
\begin{subfigure}[t]{.28\textwidth}
\includegraphics[width=\textwidth]{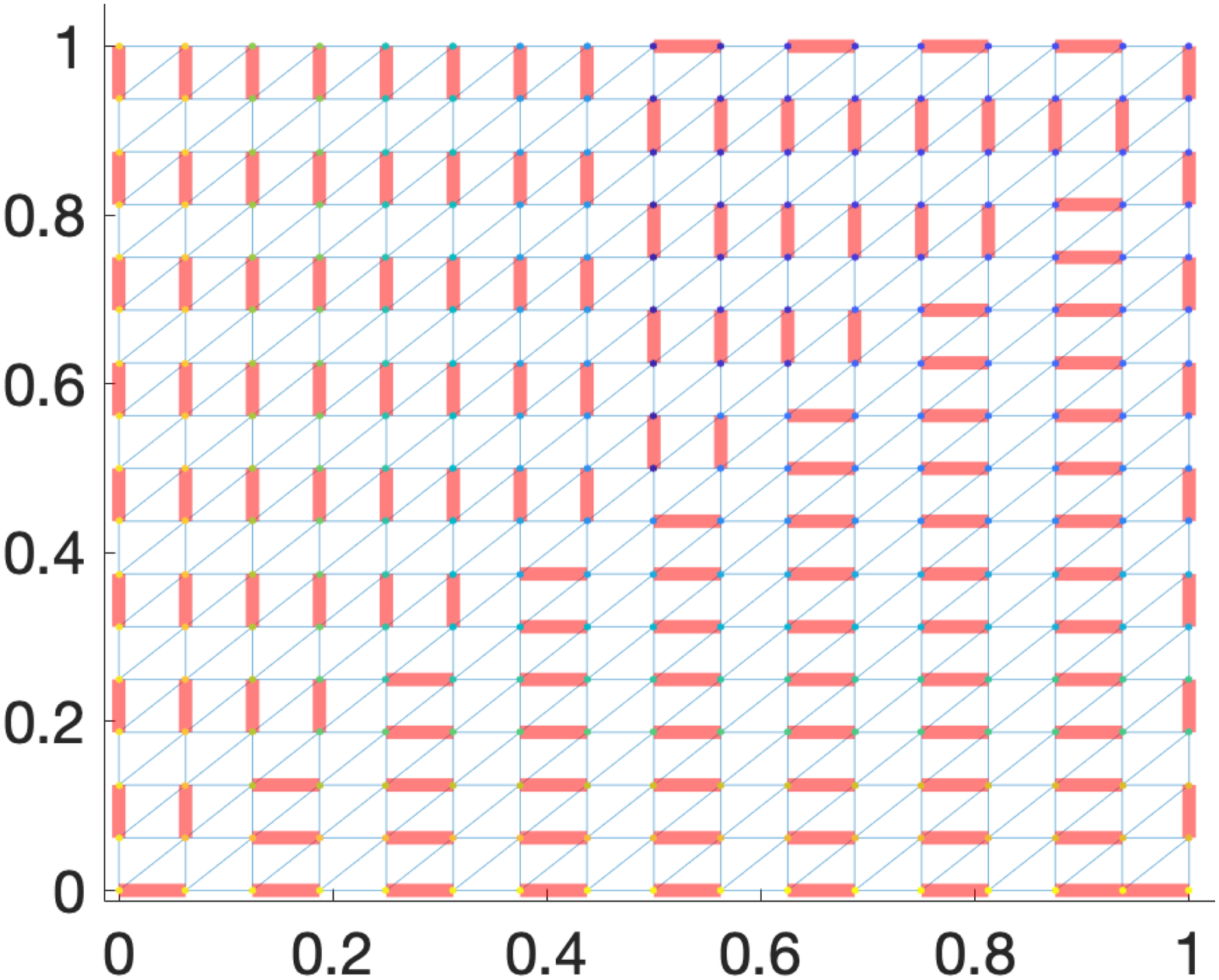}
\subcaption*{(e) Path cover aggregation based on a posteriori error estimator. }
\end{subfigure}
\begin{subfigure}[t]{.32\textwidth}
\includegraphics[width=\textwidth]{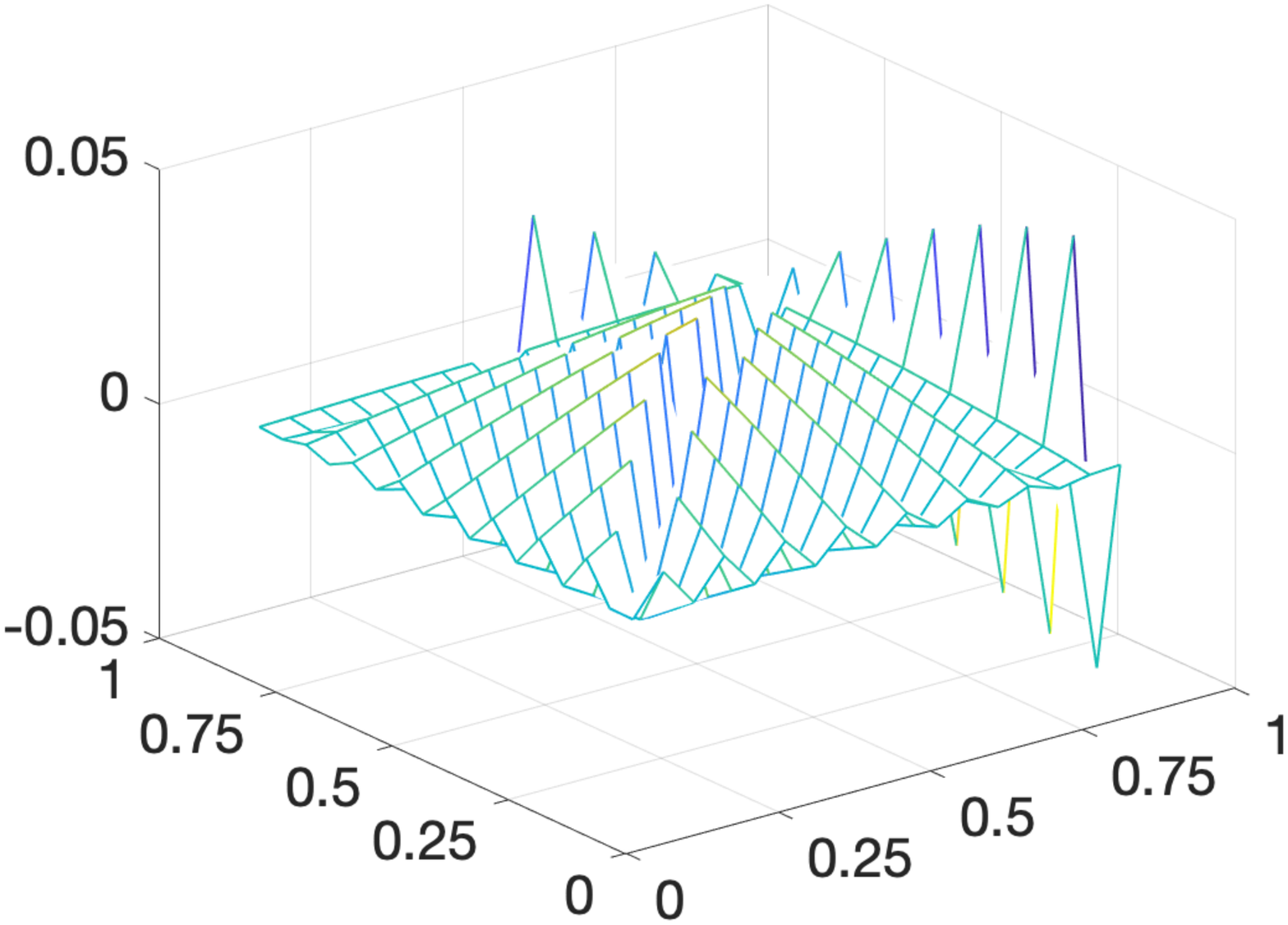}
\subcaption*{(f) Difference between the smooth error and the error after restriction and prolongation with aggregation based on error estimates.}
\end{subfigure}
\caption*{Figure 4: Path cover aggregation generated with the exact error and a posteriori error estimates, and the difference between the smooth error and the error after restriction and prolongation. The aggregation based on the error estimator assembles the one based on the exact error, but the formal one is computable at a modest cost. }
\label{fig: aggregation}
\end{figure}
 
\section{Conclusions}\label{sec:conc}
In this paper we proposed an a posteriori error estimator for solving linear systems of graph Laplacians. A novel approach is devised to reduce the computation cost of computing such an estimator in comparison to existing approaches. For sparse graphs this novel estimator can be calculated in nearly-linear time. Our approach is based on the Helmholtz decomposition on the graphs. It includes solving a linear system on a spanning tree and solving (approximately) a minimization problem in the cycle space of the graph. 

In the future, we plan to incorporate this error estimator in the adaptive AMG coarsening schemes. For example, as briefly discussed in~\cref{sec:PC-AMG}, the proposed estimator can be used in the path-cover adaptive AMG proposed in~\cite{Hu2018_PathCover}.

\bibliographystyle{siamplain}
\bibliography{references}

\end{document}